\documentclass[graybox]{svmult}

\usepackage{color}
\usepackage{amsmath}
\usepackage{amsfonts}
\usepackage{amssymb}
\usepackage{amsbsy}
\usepackage{bm}
\usepackage{graphicx}
\usepackage{extarrows}
\usepackage{tikz}
\usepackage{hyperref}

\newcommand{\formComma}{\,\text{,}}
\newcommand{\formPeriod}{\,\text{.}}
\newcommand{\R}{\mathbb{R}}%
\newcommand{\xb}{\mathbf{x}}%
\newcommand{\ub}{\mathbf{u}}%
\newcommand{\tub}{\tilde{\ub}}%
\newcommand{\tu}{\tilde{u}}%
\newcommand{\vb}{\mathbf{v}}%
\newcommand{\eb}{\mathbf{e}}%
\newcommand{\gb}{\mathbf{g}}%
\newcommand{\Xb}{\mathbf{x}}% parametrization X:{(theta,phi)} --> R3

\newcommand{\scalarprod}[1]{\big\langle{#1}\big\rangle}%
\newcommand{\exd}{\mathbf{d}} %exterior derivative
\newcommand{\lie}{\mathcal{L}} %Lie-Ableitung
\newcommand{\surf}{\mathcal{S}}
\newcommand{\gaussianCurvature}{\kappa}
\newcommand{\Grad}{\operatorname{grad}}
\newcommand{\Div}{\operatorname{div}}%
\newcommand{\Rot}{\operatorname{rot}}%
\newcommand{\DivSurf}{\Div_{\surf}}%
\newcommand{\GradSurf}{\Grad_{\surf}}
\newcommand{\RotSurf}{\Rot_{\surf}}

\newcommand{\vecLaplace}{\boldsymbol{\Delta}}
\newcommand{\laplaceDeRham}{\vecLaplace^{\textup{dR}}}
\newcommand{\laplaceRotRot}{\vecLaplace^{\textup{RR}}}
\newcommand{\laplaceGradDiv}{\vecLaplace^{\textup{GD}}}
\newcommand{\LaplaceDeRham}{Laplace-deRham }

\newcommand{\Tangent}{\mathsf{T}}
\newcommand{\vect}[1]{\mathbf{#1}}
\newcommand{\tensor}[1]{\mathbf{#1}}
\newcommand{\landau}{\mathcal{O}} %landau symbol O
\newcommand{\SC}{\mathcal{K}} % simplicial complex
\newcommand{\Vs}{\mathcal{V}} % set of vertices
\newcommand{\Es}{\mathcal{E}} % set of edges
\newcommand{\Fs}{\mathcal{T}} % set of faces
\newcommand{\face}{T} % one face
\newcommand{\FormSpace}{\Lambda^{1}} % space of 1 forms
\newcommand{\U}{u} %Komponenten des Geschwindigkeitfeldes
\newcommand{\Ub}{\mathbf{\U}} %Geschwindigkeitsfeld
\newcommand{\tU}{\tilde{u}} %Komponenten des Geschwindigkeitfeldes
\newcommand{\tUb}{\mathbf{\tU}} %Geschwindigkeitsfeld
\newcommand{\lc}{\mathbf{E}} % Levi-Civita-Tensor

\newcommand{\tensorsf}[3]{{{#1}^{#2}\!}_{#3}}

\begin{document}

\title*{Discrete exterior calculus (DEC) for the surface Navier-Stokes equation} 
\author{Ingo Nitschke, Sebastian Reuther and Axel Voigt}
\institute{
	Ingo Nitschke \at Institute of Scientific Computing, TU Dresden, 01062 Dresden, Germany, \email{ingo.nitschke@tu-dresden.de}
	\and 
	Sebastian Reuther \at Institute of Scientific Computing, TU Dresden, 01062 Dresden, Germany, \email{sebastian.reuther@tu-dresden.de}
	\and
	Axel Voigt \at Institute of Scientific Computing, TU Dresden, 01062 Dresden, Germany, Dresden Center for Computational Materials Science (DCMS), TU Dresden, 01062 Dresden, Germany and Center for Systems Biology Dresden (CSBD), Pfotenhauerstr. 108, 01307 Dresden, Germany, \email{axel.voigt@tu-dresden.de}
}

\maketitle

\abstract{
We consider a numerical approach for the incompressible surface Navier-Stokes equation. The approach is based on the covariant form and uses discrete exterior calculus (DEC) in space and a semi-implicit discretization in time. The discretization is described in detail and related to finite difference schemes on staggered grids in flat space for which we demonstrate second order convergence. We compare computational results with a vorticity-stream function approach for surfaces with genus $g(\surf) = 0$ and demonstrate the interplay between topology, geometry and flow properties. Our discretization also allows to handle harmonic vector fields, which we demonstrate on a torus.
}

\section{Introduction}

We consider a compact smooth Riemannian surface $\surf$ without boundary and an incompressible surface Navier-Stokes equation
\begin{align}
	\label{eq1}
	\partial_{t}\vb + \nabla_{\vb} \vb &= - \GradSurf p + \frac{1}{\text{Re}} \left(- \laplaceDeRham \vb + 2 \gaussianCurvature \vb \right) \\
	\label{eq2}
	\DivSurf \vb &= 0  
\end{align}
in $\surf\times\left( 0,\infty \right)$ with initial condition \( \vb \left( \xb, t=0  \right) = \vb_{0}(\xb) \in \Tangent_{\xb}\surf \). Thereby $\vb(t) \in \Tangent
\surf$ denotes the tangential surface velocity, $p(\xb,t) \in \R$ the surface pressure, $\text{Re}$ the surface Reynolds number, $\gaussianCurvature$ the Gaussian curvature,
$\Tangent_{\xb}\surf $ the tangent space on $ \xb \in \surf$, $\Tangent \surf = \cup_{\xb \in \surf} \Tangent_{\xb} \surf$ the tangent bundle and $\nabla_{\vb}, \DivSurf$ and
$\laplaceDeRham$ the covariant directional derivative, surface divergence and surface Laplace-DeRham operator, respectively. As in flat space the equation results from conservation of mass
and (tangential) linear momentum. However, differences are found in the appearing operators and the additional term including the Gaussian curvature. The Laplace-DeRham operator and the
Gaussian curvature term thereby result from the divergence of the deformation tensor and the non-commutativity of the second covariant derivative in curved spaces, see e.g.
\cite{MitreaTaylor_MA_2001,ArroyoDeSimone_PRE_2009}. The unusual sign in front of the Laplacian results from the definition of the Laplace-DeRham operator \cite{Marsden_1988}, see Section \ref{sec2}. Alternatively, the equations can also be derived from the Rayleigh dissipation potential \cite{Doerriesetal_PRE_1996}. The equations are related to the Boussinesq-Scriven constitutive law for the surface viscosity in two-phase flow problems \cite{Scriven_CES_1960,Secombetal_QJMAM_1982,Botheetal_JFM_2010} and to fluidic biomembranes  \cite{Huetal_PRE_2007,ArroyoDeSimone_PRE_2009,Fanetal_JCP_2010,Barrettetal_PRE_2015}. 
Further applications can be found in computer graphics \cite{Elcottetal_2007,Mullenetal_2009,Vaxmanetal_2016}.

While a huge literature exists for the two-dimensional Navier-Stokes equation in flat space, results for its surface counterpart eqs. \eqref{eq1} and \eqref{eq2} are rare. For treatments in
the mathematical literature we refer to \cite{EbinMarsden_AM_1970,MitreaTaylor_MA_2001}. Numerical approaches are considered in \cite{Nitschkeetal_JFM_2012,Reutheretal_MMS_2015}, where a surface vorticity-stream function formulation is introduced. 
%This follows by considering the velocity $\mathbf{v}_3 = (v_1,v_2,0)$ with the third component in direction of the surface normal $\mathbf{n}(\mathbf{x})$ and $\mathbf{v}_3 = \nabla \times \psi$ with the surface stream function $\psi$. 
This follows by considering the velocity \( \vb \) as the curl of a smooth scalar valued function \( \psi \), i.e. \( \vb = \RotSurf\psi \).
For the correct definition of the curl operator \( \RotSurf(\cdot) \) we refer to \cite{Nestleretal_arXiv_2016}.
On a compact, boundaryless, oriented Riemannian manifold of genus \( g(\surf) = 0 \), this representation is unique up to a constant by the Hodge decomposition theorem \cite{Marsden_1988}.
The resulting equations, after taking the curl and written as a system of two second order scalar surface partial differential equations, read
\begin{align} 
	\partial_t \phi + J(\psi, \phi) &= \mu (\Delta_\surf \phi + 2\DivSurf (\gaussianCurvature \GradSurf \psi)) \label{eq:vor1} \\
	\phi &= \Delta_\surf \psi \label{eq:vor2}
\end{align}
in $\surf\times\left( 0,\infty \right)$ with initial condition \( \psi \left( \xb, t=0  \right) = \psi_{0}(\xb) \in \R \). Here $\phi$ is the surface vorticity, $\Delta_\surf$ the
Laplace-Beltrami operator and $J(\psi, \phi) = \left\langle \RotSurf \psi , \GradSurf \phi \right\rangle$ the Jacobian. Eqs. \eqref{eq:vor1} and \eqref{eq:vor2} are either solved using the
surface finite element approach \cite{DziukElliott_IMAJNA_2007,DziukElliott_JCM_2007,VeyVoigt_CVS_2007}, see \cite{Nitschkeetal_JFM_2012,Reutheretal_MMS_2015} for details, or the diffuse
interface approach \cite{Raetzetal_CMS_2006}, see \cite{Reutheretal_JCP_2016} for details. The equations, but without the Gaussian curvature term, has also been discretized using a discrete exterior calculus (DEC) approach \cite{Mohamedetal_JCP_2016}. We are not aware of any direct numerical approach for eqs. \eqref{eq1} and \eqref{eq2}, which will be the purpose of this paper. Such an approach will be desirable for surfaces with genus $g(\surf) \neq 0$, as it allows to also deal with harmonic vector fields. We will introduce a DEC approach and validate the results against a surface finite element discretization for the vorticity-stream function formulation in eqs. \eqref{eq:vor1} and \eqref{eq:vor2} on surfaces with $g(\surf) = 0$ and show nontrivial solutions with $\DivSurf \vb =0$ and $\RotSurf \vb = 0$ on a torus. 

The paper is organized as follows. In Section \ref{sec2} we introduce the necessary notation and provide the formulation in covariant form. In Section \ref{sec3} the DEC discretization is described in detail and compared with known discretizations in flat space. After some analytical results for the surface Navier-Stokes equation we use the properties of Killing vector fields to validate the approach on various surfaces and demonstrate the strong interplay of geometric and vortex interactions in Section \ref{sec4}. Conclusions are drawn in Section \ref{sec5}. In the Appendix we provide additional notation and prove second order convergence for the corresponding finite difference scheme in flat space.

\section{Formulation in covariant form}
\label{sec2}

For the readers convenience we here briefly review the basic notion. A more detailed description can be found in \cite{Nestleretal_arXiv_2016}. 
The key ingredient for a covariant formulation in local coordinates \( \phi \) and \( \theta \) is the positive definite metric tensor
\begin{align}
	\tensor{g} =
		\begin{bmatrix}
			g_{\phi\phi} & g_{\phi\theta}\\
			g_{\phi\theta} & g_{\theta\theta}
		\end{bmatrix}
		= g_{\phi\phi}\,d\phi^{2} + 2g_{\phi\theta}\,d\phi \,d\theta + g_{\theta\theta}\,d\theta^{2} \formPeriod
\end{align}
\( \tensor{g}\) can be obtained from a surface parametrization $\Xb:\R^{2}\supset U \rightarrow \R^{3};\ \left( \phi, \theta \right)\mapsto \Xb \left(\phi, \theta \right)$, which maps local
coordinates to the embedded \( \R^{3} \) representation of the surface \( \surf = \Xb(U) \). The covariant components of the metric tensor are given by \( \R^{3} \) inner products of
partial derivatives of \( \Xb \), i.e. \( g_{ij} = \partial_{i}\Xb\cdot\partial_{j}\Xb \). The components of the inverse tensor \( \tensor{g}^{-1} \) are denoted by \( g^{ij} \) and the
determinant of \( \tensor{g} \) by \( \left| \tensor{g} \right| \). We denote by \( \left\{ \partial_{\phi}\Xb,\partial_{\theta}\Xb  \right\}\) the canonical basis to describe the
contravariant (tangential) vector \( \vb(\xb) \in \Tangent_{\Xb}\surf \), i.e. \( \vb(\xb) = (u^\phi,u^\theta) = u^{\phi} \partial_{\phi}\Xb + u^{\theta}\partial_{\theta}\Xb \) at a point
\( \Xb\in\surf \). Furthermore, with the arising dual basis \( \left\{ d\phi, d\theta \right\} \) we are able to write an arbitrary 1-form (covariant vector) \( \ub(\xb) \in
\Tangent^{*}_{\xb}\surf \) as \( \ub(\xb) = u_{\phi}d\phi + u_{\theta}d\theta \). This identifier choice of the covariant vector coordinates \( u_{i} \) in conjunction with representation
of \( \vb \) as above implies that \( \ub \) and \( \vb \) are related by \( \ub = \vb^{\flat} \) and \( \vb = \ub^{\sharp} \), respectively. Explicitly lowering and rising the indices can
be done using the metric tensor \( \tensor{g} \) by  \( u_{i} = g_{ij}u^{j} \) and \( u^{i} = g^{ij}u_{j} \), respectively. The scalar $p(\xb)$ is also considered as a 0-form.

We now use exterior calculus (EC) to describe all present first order differential operators by the Hodge star \( * \) and the exterior derivative \( \exd \), 
which arise algebraically (see \cite{Marsden_1988} for details).  In \cite{Marsden_1988} the \LaplaceDeRham operator is defined for \( k \)-forms on a \( n \)-dimensional Riemannian manifold by \( \laplaceDeRham:= \left( -1 \right)^{nk+1}\left( *\exd*\exd  + \exd*\exd*\right)\). For vector fields the \LaplaceDeRham operator can thus be defined canonically as composition \( (\sharp\circ\laplaceDeRham\circ\flat) \). This leads to  $\laplaceDeRham \vb = -\left(\laplaceRotRot  +  \laplaceGradDiv \right) \vb$, with the Rot-Rot-Laplace $\laplaceRotRot \vb := \RotSurf\RotSurf\vb$ and Grad-Div-Laplace $\laplaceGradDiv \vb := \GradSurf\DivSurf\vb$.  Due to the incompressibility constraint $\DivSurf \vb = 0$ we thus have $\laplaceDeRham \vb = - \laplaceRotRot \vb$ and therefor in our case only $\laplaceDeRham \ub = - (*\exd * \exd) \ub$. Eqs. \eqref{eq1} and \eqref{eq2} read in their covariant form
\begin{align}
	\label{eq3}
	\partial_{t}\ub + \nabla_{\ub}\ub &= - \exd p + \frac{1}{\text{Re}} \left( (*\exd * \exd) \ub + 2 \gaussianCurvature \ub \right) \\
	\label{eq4}
	* \exd * \ub &= 0  
\end{align}
in $\surf \times (0,\infty)$ with \( \left[ \nabla_{\ub}\ub \right]_{i} = u^{j} u_{i|j} \) to be discussed below and initial conditions $\ub(\xb,t = 0) = \ub_0(\xb) \in \Tangent^{*}_{\xb}\surf$.

\section{DEC discretization}
\label{sec3}

The mathematical foundation of discrete exterior calculus (DEC) can be found in \cite{Mercat_CMP_2001,Hirani_2003,Desbrunetal_arXiv_2005}. It follows by successively utilizing a discrete version of the Hodge star $*$ and the Stokes theorem for the exterior derivative $\exd$. The approach has been successfully used in computer graphics, e.g. surface parametrization, 
see e.g. \cite{Guetal_2003,Gortleretal_2006,Tongetal_2006}, and vector field decomposition and smoothing, see e.g. \cite{Tongetal_2003,Polthieretal_2003,Fisheretal_2006}. A rigorous treatment of the connection between discrete and continuous settings is given in \cite{Arnoldetal_AN_2006}. We discuss the discretization for each term, introduce a time discretization and compare the resulting discrete system with known discretization schemes in flat space. However, we first introduce the degrees of freedom (DOFs) and rewrite the advection term to be suitable for the DEC discretization.  

\subsection{Degrees of freedom (DOFs)}

We consider a simplicial complex \( \SC=\Vs\sqcup\Es\sqcup\Fs \) containing sets of vertices \( \Vs \), edges \( \Es \) and (triangular) faces \( \Fs \) which approximate $\surf$. The
quantities of interest in our DEC discretization are 0- and 1-forms, $p \in \Lambda^{0}(\surf)$ with $p(\xb) \in \R$ and  \( \ub \in \FormSpace(\surf) = \Tangent^{*}\surf\), 
respectively. The discrete 0-forms are considered on $v \in \Vs$, $p_h(v) := p(\xb)_{|\xb = v}$. For 1-forms we introduce DOFs as integral values on the edges \( e\in\Es \), i.e. $u_{h}(e) := \int_{\pi(e)}\ub \formComma$ with the gluing map \( \pi:\Es\rightarrow\surf \), which projects geometrically the edge \( e \) to the surface \( \surf \). The mapping \(u_{h}\in\FormSpace_{h}(\SC) \) is called the discrete 1-form of \( \ub \), since \(  u_{h}(e) \) approximates \( \ub(e)\equiv\ub(\eb) = \scalarprod{\vb,\eb} \) on an intermediate point \( \xi\in\pi(e)\subset\surf \), where the edge vector \( \eb \) exists in \( \Tangent_{\xi}\surf|_{\pi(e)} \) by the mean value theorem. Therefore, we approximate 1-forms on the restricted dual tangential space \( \Tangent_{\xi}\surf|_{\pi(e)} \), which is a one dimensional vector space in \( \xi\in\surf \) likewise the space of discrete 1-forms \( \FormSpace_{h}\left( \SC \right)|_{e} = \FormSpace_{h}\left( \{e\} \right)\) restricted to the edge \( e \), 
see \cite{Nestleretal_arXiv_2016} for details. 
Furthermore, a discrete 1-form \(  u_{h}(e) \) can be approximated as \( |e|\scalarprod{\vb,\eb}(c(e)) \) by the midpoint rule, 
with the midpoint \( c(e) = \frac{v_{1} + v_{2}}{2} \) of the edge \( e = [v_{1},v_{2}] \).
If the mesh is considered to be flat and the faces are considered to be squares, we obtain the same DOF positions as for discretizations on a staggered grid, see Section \ref{sec:staggeredGrid}.

\subsection{Approximation of the advection term}

The advection term in eq. \eqref{eq3} is not yet written in an appropriate form for a DEC discretization. We linearize this term using a Taylor expansion with a known 1-form \( \tub \) and obtain
\begin{align*}
	\left[ \nabla_{\ub}\ub  \right]_{i} = u^{j}u_{i|j} &\approx \tu^{j}\tu_{i|j} + \tu_{i|j}\left( u^{j} - \tu^{j} \right) + \tu^{j}\left( u_{i|j} - \tu_{i|j} \right) \\
	&= u^{j}\tu_{i|j} + \tu^{j}u_{i|j} - \tu^{j}\tu_{i|j}  = \left[ \nabla_{\tub}\ub + \nabla_{\ub}\tub - \nabla_{\tub}\tub\right]_{i}\formPeriod
\end{align*}
With the Levi-Cevita-tensor \( \lc \) (defined by the volumetric form \( \lc(\ub,\tub) = \mu(\ub,\tub) \)) we obtain for \( \ub,\tub \in\Tangent^* \surf \)
\begin{align*}
	\left( \RotSurf\ub \right)\left[ *\tub \right]_{i} &= E_{il}E_{jk}\tu^{l}u^{j|k} = \left( g_{ij}g_{lk} - g_{ik}g_{lj} \right)\tu^{l}u^{j|k} = \tu^{l}\left( u_{i|l} - u_{l|i} \right) \\
	&= \left[ \nabla_{\tub}\ub \right]_{i} - \left[ \tub^{\sharp}\cdot\GradSurf\ub \right]_{i}.
\end{align*}
We further have 
\begin{align*}
	\left[ \tUb^{\sharp}\cdot\GradSurf\Ub + \Ub^{\sharp}\cdot\GradSurf\tUb\right]_{i} &= \tU^{l}\U_{l|i} + \U_{l}\tensorsf{\tU}{l}{|i} = \left( \tU^{l}\U_{l} \right)_{|i} = \partial_{i}\left\langle \tUb, \Ub \right\rangle
\end{align*}
and thus also $ 2\left[ \Ub^{\sharp}\cdot\GradSurf\Ub \right]_i = \partial_{i}\left\| \Ub \right\|^{2}$. Putting everything together and using \( \RotSurf\tUb = -\DivSurf(*\tUb) \) we thus obtain
\begin{align*}
	\nabla_{\Ub}\Ub &\approx \exd\left( \left\langle \ub, \tub \right\rangle - \frac{1}{2}\left\| \tub \right\|^{2} \right) + \left( \RotSurf\ub - \RotSurf\tub \right)(*\tub) - \DivSurf(*\tub)(*\ub)
\end{align*}
which provides a suitable form for a DEC approach. By using $\RotSurf \ub =* \exd \ub$ and $\DivSurf \ub = * \exd * \ub$ we obtain
\begin{align}
	\nabla_{\Ub}\Ub &\approx \exd\left( \left\langle \ub, \tub \right\rangle - \frac{1}{2}\left\| \tub \right\|^{2} \right) + \left( * \exd \ub - * \exd \tub \right)(*\tub) - (*\exd *)(*\tub)(*\ub)
	\label{eq:Lin2}
\end{align}
which will be used for discretization.

\subsection{Time-discrete equations}
We consider a semi-implicit Euler discretization and use the approximation of the advection term with $\tub = \ub_{k}$, the solution at time $t_k$. For  \( \tau_{k} := t_{k+1} - t_{k} \)
and initial condition \( \Ub_{0} \) we get a sequence of linear systems for $k= 0,1,2,\ldots$. We introduce the generalized pressure $q_{k+1} = p_{k+1} + \left\langle \ub_{k+1},
\ub_{k} \right\rangle - \frac{1}{2}\left\| \ub_k \right\|^{2}$ and solve for $\ub_{k+1}$, $q_{k+1}$ and $p_{k+1}$
\begin{align}
	\frac{1}{\tau_{k}}\ub_{k+1} + \exd q_{k+1} + (* \exd \ub_{k+1})(* \ub_{k}) - (*\exd * )&(* \ub_{k})(* \ub_{k+1}) \nonumber \\ 
	- \frac{1}{\text{Re}} \left(  (*\exd * \exd) \ub_{k+1} + 2\gaussianCurvature \ub_{k+1} \right) &= \frac{1}{\tau_{k}}\ub_{k} + (* \exd\ub_{k})(* \ub_{k}) \label{eq:NS1DisTime}\\
	\left\langle \ub_{k+1}, \ub_{k} \right\rangle + p_{k+1} - q_{k+1} &= \frac{1}{2}\left\| \ub_{k} \right\|^{2} \\
	* \exd * \ub_{k+1} &= 0 {\color{white}\frac{1}{1}}
\end{align}
on \( \surf \).

\subsection{Fully-discrete equations}
\begin{figure}[t]
	\centering
	\begin{tikzpicture}[>=latex, line width=1pt, scale=1.2]
  % Coords
\coordinate (V0) at (2,0);
\coordinate (V1) at (2,2);
\coordinate (V2) at (0,1);
\coordinate (V3) at (4,1);
\coordinate (V4) at (4,-1);
\coordinate (V5) at (0,-1);
\coordinate (V6) at (2,-2);

%circumcenter
\coordinate (CC1) at (1.333,1);
\coordinate (CC0) at (2.666,1);
\coordinate (CC3) at (1.333,-1);
\coordinate (CC2) at (0.666,0);
\coordinate (CC4) at (2.666,-1);
\coordinate (CC5) at (3.333,0);

\coordinate (C01) at (2,1);
\coordinate (C02) at (1,0.5);
\coordinate (C03) at (3,0.5);
\coordinate (C05) at (1,-0.5);
\coordinate (C04) at (3,-0.5);
\coordinate (C06) at (2,-1);

\coordinate (C12) at (1,1.5);
\coordinate (C25) at (0,0);
\coordinate (C56) at (1,-1.5);
\coordinate (C46) at (3,-1.5);
\coordinate (C34) at (4,0);
\coordinate (C13) at (3,1.5);

\fill[black, opacity=0.2] (CC0) -- (CC1) -- (CC2) -- (CC3) -- (CC4) -- (CC5) -- (CC0);
\fill[black, opacity=0.3] (V0) -- (CC5) -- (CC0) -- (V0);
%\fill[opacity=0.4]  (V1) -- (V3) -- (V4) -- (V6) -- (V5) -- (V2);
  % Arrows\tilde{\sigma}
\draw[](V0) -- (V1);
\draw[](V1) --  (V2);
\draw[](V0) -- (V2);
\draw[](V1) --  (V3);
\draw[blue,->](V0) -- (V3);
\draw[](V0) -- (V4);
\draw[](V0) -- (V5);
\draw[](V5) --  (V2);
\draw[](V3) --  (V4);
\draw[](V4) -- (V6);
\draw[](V0) -- (V6);
\draw[](V6) -- (V5);
 
\draw (CC0) -- (CC1) -- (CC2) -- (CC3) -- (CC4) -- (CC5) -- (CC0);
\draw[dotted] (C12) -- (CC1);
\draw[dotted] (C25) -- (CC2);
\draw[dotted] (C56) -- (CC3);
\draw[dotted] (C46) -- (CC4);
\draw[dotted] (C34) -- (CC5);
\draw[dotted] (C13) -- (CC0);
\draw[dotted] (CC1) -- (CC4);
\draw[dotted] (CC2) -- (CC5);
\draw[dotted] (CC3) -- (CC0);

\draw[dotted] (CC1) -- (V1);
\draw[dotted] (CC2) -- (V2);
\draw[dotted] (CC3) -- (V5);
\draw[dotted] (CC4) -- (V6);
\draw[dotted] (CC5) -- (V4);
\draw[dotted] (CC0) -- (V3);

\draw[dotted] (CC1) -- (V2);
\draw[dotted] (CC2) -- (V5);
\draw[dotted] (CC3) -- (V6);
\draw[dotted] (CC4) -- (V4);
\draw[dotted] (CC5) -- (V3);
\draw[dotted] (CC0) -- (V1);

%\draw[blue,->] (CC5) -- (CC0); 

\fill[red] (V0) node[left] {\(v\ \)} circle (2pt);
\fill (V1) circle (2pt);
\fill (V2) circle (2pt);
\fill (V3) circle (2pt);
\fill (V4) circle (2pt);
\fill (V5)circle (2pt);
\fill (V6)circle (2pt);

\fill (CC0) circle (2pt);
\fill (CC1) circle (2pt);
\fill (CC2) circle (2pt);
\fill (CC3) circle (2pt);
\fill (CC4) circle (2pt);
\fill (CC5) circle (2pt);

\fill (C01) circle (2pt);
\fill (C02) circle (2pt);
\fill (C03) node[right, blue] {\(\ e\)}circle (2pt);
\fill (C05) circle (2pt);
\fill (C04) circle (2pt);
\fill (C06) circle (2pt);

\fill (C12) circle (2pt);
\fill (C25) circle (2pt);
\fill (C56) circle (2pt);
\fill (C46) circle (2pt);
\fill (C34) circle (2pt);
\fill (C13) circle (2pt);

\end{tikzpicture}
	\caption{Circumcentric subdivision of a simple simplicial complex around a vertex \( v \). The Voronoi cell \( \star v \) is marked light gray and the part \( A_{ve} \) dark gray.}
	\label{fig:VCellDecomp}
\end{figure}
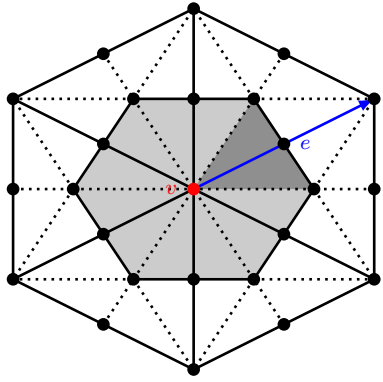

 The used notation follows \cite{Nestleretal_arXiv_2016}, see also Appendix \ref{app:a}. For the discrete 0-forms $p_h, q_h \in \Lambda_h^0(\SC)$, 1-forms \( u_{h}\in\FormSpace_{h}(\SC) \), sign mappings \( s_{\circ,\circ}\in\left\{ -1,+1 \right\} \), volumes \( \left| \cdot \right| \), Voronoi cells \( \star v \), Voronoi edges \( \star e \) and the ``belongs-to'' relations \( \succ \) and \( \prec\) we obtain for $\laplaceRotRot$, $\DivSurf$ and $\RotSurf$
\begin{align*}
	 ( *\exd * \exd) \ub (e) &\approx -\frac{\left| e \right|}{\left| \star e \right|} \sum_{\face\succ e} \frac{s_{\face,e}}{\left| \face \right|} \sum_{\tilde{e}\prec\face} s_{\face,\tilde{e}}\, u_{h}(\tilde{e}) \formComma\\
	(*\exd * \ub)(v) &\approx -\frac{1}{|\star v|}\sum_{\tilde{e}\succ v}s_{v,\tilde{e}}\frac{|\star \tilde{e}|}{|\tilde{e}|}\ub(\tilde{e})\formComma \\
	(*\exd\ub) (c(e)) &  \approx  \frac{1}{\sum_{\face\succ e}|\face|}\sum_{\face\succ e}(\exd\ub)(\face) = \frac{1}{\sum_{\face\succ e}|\face|}\sum_{\face\succ
  e}\sum_{\tilde{e}\prec\face}s_{\face,\tilde{e}}\ub(\tilde{e})
\end{align*}
respectively. 
The last line above follows from a special Hodge dualism between midpoint \( c(e) \) and face union \( \bigcup_{\face\succ e}\face =: \hat{\star}c(e) \), 
such that function evaluations at \( c(e) \) are integral mean values over \(\hat{\star}c(e)  \).
This allows to approximate
\begin{align*}
	(* \exd \ub)(*\tub)(e) &\approx \frac{(*\tilde{\ub})(e)}{\sum_{\face\succ e}|\face|}\sum_{\face\succ e}\sum_{\tilde{e}\prec\face}s_{\face,\tilde{e}}\ub(\tilde{e}) \formComma \\
	(*\exd * \tub)(* \ub)(e) &\approx -\frac{1}{2}\left( \sum_{v\prec e} \frac{1}{|\star v|}\sum_{\tilde{e}\succ v}s_{v,\tilde{e}}\frac{|\star \tilde{e}|}{|\tilde{e}|}\tilde{\ub}(\tilde{e}) \right) (*\ub)(e) \formPeriod
\end{align*}
With the Stokes theorem we further have $(\exd q)(e) = q(v_j) - q(v_i)$ for $e = [v_i,v_j]$. What remains to define is a discrete Hodge operator and a discrete version of the inner product. We approximate
\begin{align*}
	\left( *\ub \right)(e) &\approx \circledast u_h(e) \\
	&= \frac{1}{4}\sum_{\face \succ e}\sum_{\substack{\tilde{e}\prec\face\\\tilde{e}\ne e}} \frac{s_{e\tilde{e}}}{\sqrt{\left| e \right|^{2}\left| \tilde{e} \right|^{2} - \left( \vect{e}\cdot\tilde{\vect{e}}\right)^{2}}} \left( \left( \vect{e}\cdot\tilde{\vect{e}} \right) u_h(e) -\left| e \right|^{2} u_h(\tilde{e}) \right) \formPeriod
\end{align*}
For other possibilities we refer to \cite{Mohamed_CAD_2016}. For the inner product we follow \cite{Nestleretal_arXiv_2016} and define
\begin{align*}
	\left\langle \tilde{\ub} , \ub \right\rangle(c(e)) &\approx \frac{1}{\left| e \right|^{2}}\left( \tilde{\ub}(e) \ub(e) + \left( *\tilde{\ub} \right)(e) \left( *\ub \right)(e) \right) \formPeriod
\end{align*}
In order to approximate the inner product at a primal vertex we consider the decomposition of the Voronoi cell $\star v = \sum_{e\succ v} A_{ve}$, with $\left|A_{ve}\right| = (\left| e \right|\left| \star e \right|)/4$, see Fig. \ref{fig:VCellDecomp}. We thus obtain
\begin{align*}
	\left\langle \tilde{\Ub} , \Ub \right\rangle(v) &\approx \frac{1}{\left| \star v \right|}\left( *\left\langle \tilde{\Ub} , \Ub \right\rangle  \right)(\star v) = \frac{1}{\left| \star v \right|}\sum_{e\succ v}\left( *\left\langle \tilde{\Ub} , \Ub \right\rangle  \right)(A_{ve}) \\
	&\approx\frac{1}{\left| \star v \right|}\sum_{e\succ v} \frac{\left| e \right|\left| \star e \right|}{4}  \left\langle \tilde{\Ub} , \Ub \right\rangle(c(e)) \\
	&\approx\frac{1}{4\left| \star v \right|}\sum_{e\succ v} \frac{\left| \star e \right|}{\left| e \right|} \left( \tilde{\Ub}(e) \Ub(e) + \left( *\tilde{\Ub} \right)(e) \left( *\Ub \right)(e) \right)
  \formPeriod
\end{align*}
In case $\kappa$ is not given analytically, a numerical approximation is required, which can effectively be done using a DEC approach for the Weingarten map \cite{Nitschkeetal_inprep}.

\subsection{Linear system}
Putting everything together and using an additional equation  
\begin{align}
	\circledast\ub_{k+1}(e) - (* \ub)_{k+1}(e) = 0
\end{align}
for all \( e\in\Es \) to determine the Hodge dual 1-form 
defines a linear system for $\ub_{k+1}, (* \ub)_{k+1} \in \Lambda_h^1(\SC)$ and $q_{k+1}, p_{k+1} \in \Lambda_h^0(\SC)$. An appropriate assembly over $e \in \Es$ and $v \in \Vs$ results in a sparse matrix $M_{k+1}\in\R^{2(|\Es| + |\Vs|)\times 2(|\Es| + |\Vs|)}$ and the right hand side vector \( r_{k}\in\R^{2(|\Es| + |\Vs|)} \). To determine the pressure we replace a row in $M_{k+1}$ and $r_{k+1}$  to ensure $p_{k+1} (v_0) = 0$ at $v_0 \in \Vs$. The linear system is solved using \( \texttt{umfpack} \).

\subsection{Comparison with finite difference schemes on uniform rectangular meshes in two dimensions} \label{sec:staggeredGrid}

To compare the resulting scheme with known discretization schemes we consider the two-dimensional Navier-Stokes equation in flat space. The Gaussian curvature $\gaussianCurvature$ vanishes and the surface operators reduce to the classical two-dimensional operators $\Grad$, $\Div$, $\Rot$ and $\Delta$. Instead of the simplicial complex \( \SC \) we consider for simplicity a uniform rectangular mesh. The DEC discretization can then be considered as introduced above.
\begin{figure}[t]
	\centering
	\begin{tikzpicture}[>=latex, line width=2pt, scale=2.5]
\usetikzlibrary{calc}

\edef\n{3}

\pgfmathparse{\n-1}
\edef\nmone{\pgfmathresult}

\pgfmathparse{\n-2}
\edef\nmtwo{\pgfmathresult}

%all the gray (dual) stuff
\begin{scope}[gray,  line width=1pt]
%dual vertices
\foreach \x in {0,...,\nmone} {
	\foreach \y in {0,...,\nmone} {
		\coordinate (V) at (\x+0.5,\y+0.5);
		\fill (V) circle(0.5pt);
	}
}

%dual x-edges
\foreach \x in {0,...,\nmtwo} {
	\foreach \y in {0,...,\nmone} {
		\coordinate (V1) at (\x+0.5,\y+0.5);
		\coordinate (V2) at (\x+1.5,\y+0.5); 
		\draw[<-] (V1) -- (V2);
	}
}

%dual y-edges
\foreach \x in {0,...,\nmone} {
	\foreach \y in {0,...,\nmtwo} {
		\coordinate (V1) at (\x+0.5,\y+0.5);
		\coordinate (V2) at (\x+0.5,\y+1.5); 
		\draw[->] (V1) -- (V2);
	}
}

%all remaining (half) dual edges on the boundaries 
\foreach \xy in {0,...,\nmone} {
	\draw[->,densely dotted] (\xy+0.5, 0) --  (\xy+0.5, 0.5); 
	\draw[->,densely dotted] (\xy+0.5, \n-0.5) --  (\xy+0.5, \n);
	\draw[<-,densely dotted] (0,\xy+0.5) -- (0.5,\xy+0.5);
	\draw[<-,densely dotted] (\n-0.5,\xy+0.5) -- (\n,\xy+0.5);
}
\end{scope}

%primal vertices
\foreach \x in {0,...,\n} {
	\foreach \y in {0,...,\n} {
		\coordinate (V) at (\x,\y);
		\fill (V) circle(1pt);
	}
}

%primal x-edges
\foreach \x in {0,...,\nmone} {
	\foreach \y in {0,...,\n} {
		\coordinate (V1) at (\x,\y);
		\coordinate (V2) at (\x+1,\y); 
		\draw[->] (V1) -- (V2);
	}
}

%primal y-edges
\foreach \x in {0,...,\n} {
	\foreach \y in {0,...,\nmone} {
		\coordinate (V1) at (\x,\y);
		\coordinate (V2) at (\x,\y+1); 
		\draw[->] (V1) -- (V2);
	}
}

%primal edge circumcenters
\foreach \xy in {0,...,\n} {
	\foreach \yx in {0,...,\nmone} {
		\coordinate (VX) at (\yx+0.5,\xy);
		\coordinate (VY) at (\xy,\yx+0.5);
		\fill (VX) circle(1pt);
		\fill (VY) circle(1pt);
	}
}

%nodes
\node[below right] at (1,1) {$v_{i,j}$};
\node[below right] at (2,1) {$v_{i+1,j}$};

\node[above right] at (1.5,1) {$e^x_{i,j}$};
\node[above right] at (1.5,0) {$e^x_{i,j-1}$};
\node[above right] at (1.5,2) {$e^x_{i,j+1}$};

\node[below right] at (1,1.5) {$e^y_{i,j}$};
\node[below right] at (2,1.5) {$e^y_{i+1,j}$};
\node[below right] at (1,0.5) {$e^y_{i,j-1}$};
\node[below right] at (2,0.5) {$e^y_{i+1,j-1}$};

\end{tikzpicture}
	\caption{Staggered grid with dual mesh and orientation. The components of the velocity $u^{x}$ and $u^{y}$ are defined on the midpoints of the edges, $u^{x}$ on the horizontal $e^x$ and
	$u^{y}$ on the vertical $e^y$, and the pressure is defined on the vertices $v$. Such meshes are also known as Arakawa C-meshes \cite{Arakawa_MCP_1977}.}
	\label{fig:staggeredGrid}
\end{figure}
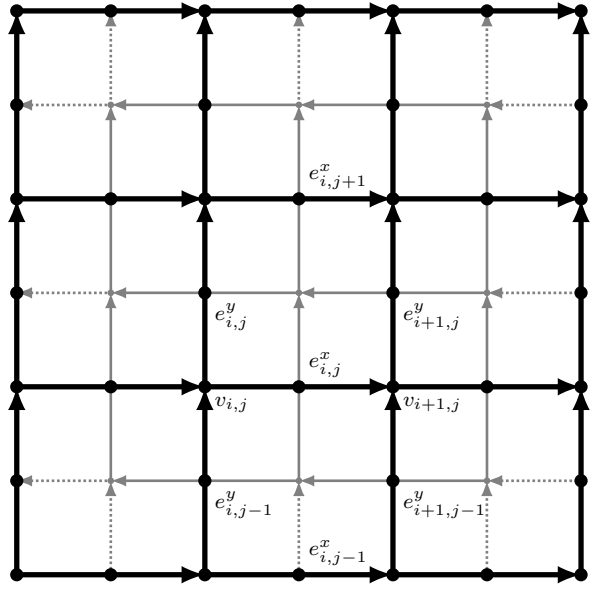

We identify the vector-components in the midpoints of the edges as the discrete 1-form $u_h$. We thus obtain with the grid spacing $h$ and the notation in Fig. \ref{fig:staggeredGrid} 
\begin{align*}
	u^{x}_{ij} := u^{x}(c(e^{x}_{i,j}))
		= \frac{1}{h} \U_{h}(e^{x}_{i,j}), &\qquad
	u^{y}_{ij} := u^{y}(c(e^{y}_{i,j}))
		= \frac{1}{h} \U_{h}(e^{y}_{i,j}) 
\end{align*}
For the pressure we obtain with the discrete 0-form $q_h$
\begin{align*}
	q_{i,j} := q(v_{i,j}) = q_{h}(v_{i,j}) \formPeriod
\end{align*}

\begin{figure}[t]
	\centering
	  \begin{tikzpicture}[>=latex, line width=1.5pt, scale=1.7]
\usetikzlibrary{calc}
\bf
\edef\n{3}

\pgfmathparse{\n-1}
\edef\nmone{\pgfmathresult}

\pgfmathparse{\n-2}
\edef\nmtwo{\pgfmathresult}

%all the gray (dual) stuff
\begin{scope}[gray,  line width=1pt]
%dual vertices
\foreach \x in {0,...,\nmone} {
	\foreach \y in {0,...,\nmone} {
		\coordinate (V) at (\x+0.5,\y+0.5);
		\fill (V) circle(0.5pt);
	}
}

%dual x-edges
\foreach \x in {0,...,\nmtwo} {
	\foreach \y in {0,...,\nmone} {
		\coordinate (V1) at (\x+0.5,\y+0.5);
		\coordinate (V2) at (\x+1.5,\y+0.5); 
		\draw[<-] (V1) -- (V2);
	}
}

%dual y-edges
\foreach \x in {0,...,\nmone} {
	\foreach \y in {0,...,\nmtwo} {
		\coordinate (V1) at (\x+0.5,\y+0.5);
		\coordinate (V2) at (\x+0.5,\y+1.5); 
		\draw[->] (V1) -- (V2);
	}
}

%all remaining (half) dual edges on the boundaries 
\foreach \xy in {0,...,\nmone} {
	\draw[->,densely dotted] (\xy+0.5, 0) --  (\xy+0.5, 0.5); 
	\draw[->,densely dotted] (\xy+0.5, \n-0.5) --  (\xy+0.5, \n);
	\draw[<-,densely dotted] (0,\xy+0.5) -- (0.5,\xy+0.5);
	\draw[<-,densely dotted] (\n-0.5,\xy+0.5) -- (\n,\xy+0.5);
}
\end{scope}

%primal vertices
\foreach \x in {0,...,\n} {
	\foreach \y in {0,...,\n} {
		\coordinate (V) at (\x,\y);
		\fill (V) circle(1pt);
	}
}

%primal x-edges
\foreach \x in {0,...,\nmone} {
	\foreach \y in {0,...,\n} {
		\coordinate (V1) at (\x,\y);
		\coordinate (V2) at (\x+1,\y); 
		\draw[->] (V1) -- (V2);
	}
}

%primal y-edges
\foreach \x in {0,...,\n} {
	\foreach \y in {0,...,\nmone} {
		\coordinate (V1) at (\x,\y);
		\coordinate (V2) at (\x,\y+1); 
		\draw[->] (V1) -- (V2);
	}
}

%primal edge circumcenters
\foreach \xy in {0,...,\n} {
	\foreach \yx in {0,...,\nmone} {
		\coordinate (VX) at (\yx+0.5,\xy);
		\coordinate (VY) at (\xy,\yx+0.5);
		\fill (VX) circle(1pt);
		\fill (VY) circle(1pt);
	}
}

%schemata
\begin{scope}[line width=3.5pt]
\edef\sh{0.6pt}
\draw[blue, yshift=\sh, xshift=\sh] (0,1.5) circle(1pt) node[above right]{+1} -- (1,1.5) circle(1pt) node[above right]{-2}  -- (2,1.5) circle(1pt)node[above right]{+1};
\draw[blue, yshift=\sh, xshift=\sh] (0.5,1)circle(1pt) node[above right]{+1} -- (0.5,2) circle(1pt)node[above right]{-1};
\draw[blue, yshift=\sh, xshift=\sh] (1.5,1)circle(1pt) node[above right]{-1} -- (1.5,2) circle(1pt)node[above right]{+1};

\draw[red] (1.5,0) circle(1pt) node[above left]{+1} -- (1.5,1)
circle(1pt) node[below left]{-2}  -- (1.5,2) circle(1pt)node[above left]{+1};
\draw[red] (1,1.5) circle(1pt) node[above left]{+1} -- (2,1.5) circle(1pt)node[above left]{-1};
\draw[red] (1,0.5) circle(1pt) node[above left]{-1} -- (2,0.5) circle(1pt)node[above left]{+1};
\end{scope}

%nodes
%\node[below right] at (1,1) {$v_{i,j}$};
%\node[below right] at (2,1) {$v_{i+1,j}$};

\node[below right] at (1.5,1) {$u^x_{i,j}$};
%\node[above right] at (1.5,0) {$e^x_{i,j-1}$};
%\node[above right] at (1.5,2) {$e^x_{i,j+1}$};

\node[below right] at (1,1.5) {$u^y_{i,j}$};
%\node[below right] at (2,1.5) {$e^y_{i+1,j}$};
%\node[below right] at (1,0.5) {$e^y_{i,j-1}$};
%\node[below right] at (2,0.5) {$e^y_{i+1,j-1}$};

\end{tikzpicture}
	  \begin{tikzpicture}[>=latex, line width=1.5pt, scale=1.7]
\usetikzlibrary{calc}

\bf
\edef\n{3}

\pgfmathparse{\n-1}
\edef\nmone{\pgfmathresult}

\pgfmathparse{\n-2}
\edef\nmtwo{\pgfmathresult}

%all the gray (dual) stuff
\begin{scope}[gray,  line width=1pt]
%dual vertices
\foreach \x in {0,...,\nmone} {
	\foreach \y in {0,...,\nmone} {
		\coordinate (V) at (\x+0.5,\y+0.5);
		\fill (V) circle(0.5pt);
	}
}

%dual x-edges
\foreach \x in {0,...,\nmtwo} {
	\foreach \y in {0,...,\nmone} {
		\coordinate (V1) at (\x+0.5,\y+0.5);
		\coordinate (V2) at (\x+1.5,\y+0.5); 
		\draw[<-] (V1) -- (V2);
	}
}

%dual y-edges
\foreach \x in {0,...,\nmone} {
	\foreach \y in {0,...,\nmtwo} {
		\coordinate (V1) at (\x+0.5,\y+0.5);
		\coordinate (V2) at (\x+0.5,\y+1.5); 
		\draw[->] (V1) -- (V2);
	}
}

%all remaining (half) dual edges on the boundaries 
\foreach \xy in {0,...,\nmone} {
	\draw[->,densely dotted] (\xy+0.5, 0) --  (\xy+0.5, 0.5); 
	\draw[->,densely dotted] (\xy+0.5, \n-0.5) --  (\xy+0.5, \n);
	\draw[<-,densely dotted] (0,\xy+0.5) -- (0.5,\xy+0.5);
	\draw[<-,densely dotted] (\n-0.5,\xy+0.5) -- (\n,\xy+0.5);
}
\end{scope}

%primal vertices
\foreach \x in {0,...,\n} {
	\foreach \y in {0,...,\n} {
		\coordinate (V) at (\x,\y);
		\fill (V) circle(1pt);
	}
}

%primal x-edges
\foreach \x in {0,...,\nmone} {
	\foreach \y in {0,...,\n} {
		\coordinate (V1) at (\x,\y);
		\coordinate (V2) at (\x+1,\y); 
		\draw[->] (V1) -- (V2);
	}
}

%primal y-edges
\foreach \x in {0,...,\n} {
	\foreach \y in {0,...,\nmone} {
		\coordinate (V1) at (\x,\y);
		\coordinate (V2) at (\x,\y+1); 
		\draw[->] (V1) -- (V2);
	}
}

%primal edge circumcenters
\foreach \xy in {0,...,\n} {
	\foreach \yx in {0,...,\nmone} {
		\coordinate (VX) at (\yx+0.5,\xy);
		\coordinate (VY) at (\xy,\yx+0.5);
		\fill (VX) circle(1pt);
		\fill (VY) circle(1pt);
	}
}

%schemata
\begin{scope}[line width=3.5pt]
\draw[blue] (0,1.5) circle(1pt) node[above left]{+1} -- (1,1.5) circle(1pt) node[above left]{-4}  -- (2,1.5) circle(1pt)node[above left]{+1};
\draw[blue] (1,0.5) circle(1pt) node[above left]{+1} -- (1,2.5) circle(1pt)node[above left]{+1};

\draw[red] (0.5,1) circle(1pt) node[above left]{+1} -- (1.5,1)
circle(1pt) node[below left]{-4}  -- (2.5,1) circle(1pt)node[above left]{+1};
\draw[red] (1.5,0) circle(1pt) node[above left]{+1} -- (1.5,2) circle(1pt)node[above left]{+1};
\end{scope}

%nodes
%\node[below right] at (1,1) {$v_{i,j}$};
%\node[below right] at (2,1) {$v_{i+1,j}$};

\node[above right] at (1.5,1) {$u^x_{i,j}$};
%\node[above right] at (1.5,0) {$e^x_{i,j-1}$};
%\node[above right] at (1.5,2) {$e^x_{i,j+1}$};

\node[below right] at (1,1.5) {$u^y_{i,j}$};
%\node[below right] at (2,1.5) {$e^y_{i+1,j}$};
%\node[below right] at (1,0.5) {$e^y_{i,j-1}$};
%\node[below right] at (2,0.5) {$e^y_{i+1,j-1}$};

\end{tikzpicture}
  \caption{(left) Illustration of the difference schemes for \( \left( \laplaceRotRot u \right)^{x} \) (red) and \( \left( \laplaceRotRot u \right)^{y} \) (blue).
           (right) Illustration of the schemes for \( \left( \vecLaplace u \right)^{x} \) (red) and \( \left( \vecLaplace u \right)^{y} \) (blue),
            which is the well known five-point stencil.}
	\label{fig:laplacesStaggeredGrid}
\end{figure}
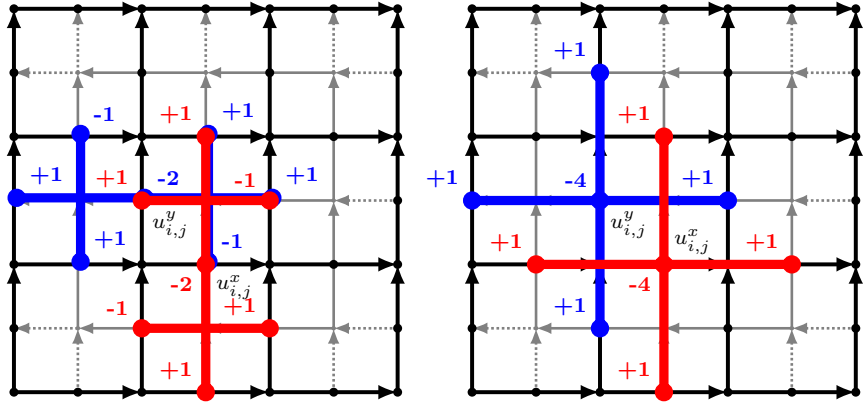

To analyze the scheme we here only consider the discretization of the Laplace operator, which is restricted to $\vecLaplace = -\laplaceRotRot$ in the present case, with
\begin{align}\label{eq:lrotrotflat}
	\left( \laplaceRotRot u \right)^{x} = \partial_{y}^{2}u^{x} - \partial_{x}\partial_{y}u^{y}, & \qquad
	\left( \laplaceRotRot u \right)^{y} = \partial_{x}^{2}u^{y} - \partial_{x}\partial_{y}u^{x}
\end{align}
and with our DEC discretization
\begin{align}\label{eq:lrotrotstaggered}
  \begin{aligned}
	&(\laplaceRotRot u)^{x}_{i,j} = \\
	&\quad\frac{1}{h^{2}}\left( u^{x}_{i,j+1} +  u^{x}_{i,j-1} - 2u^{x}_{i,j}
						+u^{y}_{i,j} - u^{y}_{i+1,j}
						+u^{y}_{i+1,j-1} -u^{y}_{i,j-1} \right) \\
	&(\laplaceRotRot u)^{y}_{i,j} = \\
	&\quad\frac{1}{h^{2}}\left( u^{y}_{i+1,j} +  u^{y}_{i-1,j} - 2u^{y}_{i,j}
						-u^{x}_{i,j} + u^{x}_{i,j+1}
						-u^{x}_{i-1,j+1} +u^{x}_{i-1,j} \right).
  \end{aligned}
\end{align}
This unusual stencil is visualized in Fig. \ref{fig:laplacesStaggeredGrid}. For the full Laplace operator $\vecLaplace = -\left(\laplaceRotRot  +  \laplaceGradDiv \right)$, as considered in
\cite{Nestleretal_arXiv_2016} and also typically used in flat space, we obtain 
\begin{align*}
	(\vecLaplace u)^{\{x,y\}}_{i,j} &= \frac{1}{h^{2}}\left( u^{\{x,y\}}_{i+1,j} +  u^{\{x,y\}}_{i-1,j} + u^{\{x,y\}}_{i,j+1} +  u^{\{x,y\}}_{i,j-1} - 4u^{\{x,y\}}_{i,j} \right)
\end{align*}
which is the usual five-point stencil, again visualized in Fig. \ref{fig:laplacesStaggeredGrid}. 
We thus have $(\vecLaplace u)^{\{x,y\}}_{i,j} \neq  (\laplaceRotRot u)^{\{x,y\}}_{i,j}$, even if the identity holds in the continuous case under the incompressibility constraint.
However, the order of consistency is \( \landau(h^{2}) \) for both stencils, which can be shown by a Taylor expansion for each component, see Appendix \ref{app:b} for details. 

The lower order terms can be compared in a similar way and lead to typical finite difference discretizations. However, a comparison of the full model strongly depends on the approximation of the advection term and will thus not be done. We conclude that the proposed DEC discretization, if considered on a uniform rectangular mesh in flat space, can be related to finite difference schemes with the same order of consistency as established approaches. Similar comparisons with finite difference schemes have also been considered for scalar valued problems in \cite{Griebeletal_SPP_2017}.
 
\section{Results}\label{sec4}

\subsection{Energy dissipation}

As in flat space we can show that $\dot{E} = \frac{d}{dt} \frac{1}{2} \int_{\surf} \left\| \ub \right\|^{2} \mu = \int_{\surf} \left\langle \Ub , \dot{\Ub} \right\rangle \mu \leq 0$. The only term which requires a remark is the viscous part $\frac{1}{\text{Re}} \left( (*\exd * \exd) \ub + 2 \gaussianCurvature \ub \right)$. 
By using the Frobenius inner product for tensors we obtain  
\begin{align*}
	&\int_{\surf} \left\langle \ub , \frac{1}{\text{Re}} \left( (* \exd * \exd) \ub + 2 \gaussianCurvature \ub \right) \right\rangle \mu = \frac{1}{\text{Re}} \int_{\surf} \left\langle \ub , \DivSurf \lie_{\ub} \gb \right\rangle \mu \\ 
	&= \frac{1}{\text{Re}}  \int_{\surf} u^{i} \left[\lie_{\ub^{\sharp}}\gb\right]_{ij}{\color{white} }^{|j} \mu 
	= -\frac{1}{\text{Re}}  \int_{\surf} u^{i|j}\left( u_{i|j} + u_{j|i} \right) \mu\\
	&= -\frac{1}{\text{Re}}  \int_{\surf} \left< \GradSurf \ub , \GradSurf \ub + (\GradSurf \ub)^{T} \right>\mu \\
	&\overset{(*)}{=} -\frac{1}{2 \text{Re}}  \int_{\surf} \left\| \GradSurf \ub + (\GradSurf \ub)^{T}  \right\|^{2} \mu
	= - \frac{1}{2 \text{Re}}  \int_{\surf} \left\| \lie_{\ub^{\sharp}}\gb  \right\|^{2} \mu \leq 0,
\end{align*}
with the Lie-derivative $\lie_{\ub^{\sharp}}$ and \( (*) \) following from the component wise computation
\begin{align*}
	& u^{i|j}\left( u_{i|j} + u_{j|i} \right) \\
	&= \frac{1}{2}\left( u^{i|j} + u^{j|i} + u^{i|j} -  u^{j|i}\right)\left( u_{i|j} + u_{j|i} \right) \\
	&= \frac{1}{2}\left[\left( u^{i|j} + u^{j|i} \right)\left( u_{i|j} + u_{j|i} \right) + \left(  u^{i|j}u_{i|j} - u^{j|i}u_{j|i} \right) + \left(  u^{i|j}u_{j|i} -  u^{j|i}u_{i|j}\right)\right] \\
	&= \frac{1}{2}\left( u^{i|j} + u^{j|i} \right)\left( u_{i|j} + u_{j|i} \right) \formPeriod
\end{align*}
As in flat space we obtain a non-dissipative system for the corresponding surface Euler equation ($\text{Re} \to \infty$). However, the system is also non-dissipative for
$\lie_{\ub^{\sharp}}\gb = 0$, so called Killing vector fields \cite{Marsden_1988}, which can be realized on rotational symmetric surfaces. We will use this property in various examples.
 
\subsection{Numerical dissipation}

\begin{figure}[t]
	\begin{minipage}{0.3\textwidth}
		\includegraphics[width=\textwidth]{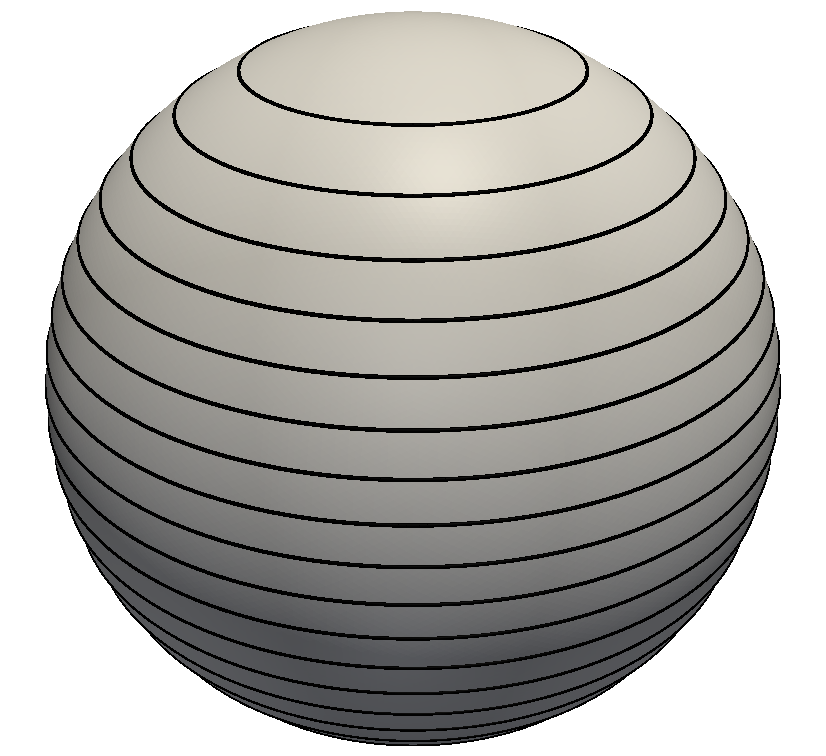}
	\end{minipage}
	\hfill
	\begin{minipage}{0.65\textwidth}
		\centering
		\begin{tabular}{lp{0.04\textwidth}p{0.24\textwidth}lp{0.04\textwidth}p{0.24\textwidth}l}
			& & \multicolumn{2}{c}{SFEM} & & \multicolumn{2}{c}{DEC} \\ 
			\hline\noalign{\smallskip}
			$h$ & & $|E_0 - E(10)|$ & $EOC$ & &  $|E_0 - E(10)|$ & $EOC$ \\
			\noalign{\smallskip}\svhline\noalign{\smallskip}
			$0.109$ & & $0.0777$ & $-$ & & $0.0737$ & $-$ \\ 
			$0.066$ & & $0.0295$ & $1.929$ & & $0.0280$ & $1.931$ \\ 
			$0.048$ & & $0.0154$ & $2.050$ & & $0.0146$ & $2.052$ \\ 
			$0.031$ & & $0.0063$ & $2.031$ & & $0.0060$ & $2.032$ \\ 
			$0.022$ & & $0.0031$ & $2.039$ & & $0.0030$ & $2.039$ \\ 
			$0.015$ & & $0.0016$ & $1.777$ & & $0.0015$ & $1.785$ \\ 
			\noalign{\smallskip}\hline\noalign{\smallskip}
		\end{tabular}
	\end{minipage} \\
	\begin{minipage}{0.49\textwidth}
		\includegraphics[width=\textwidth]{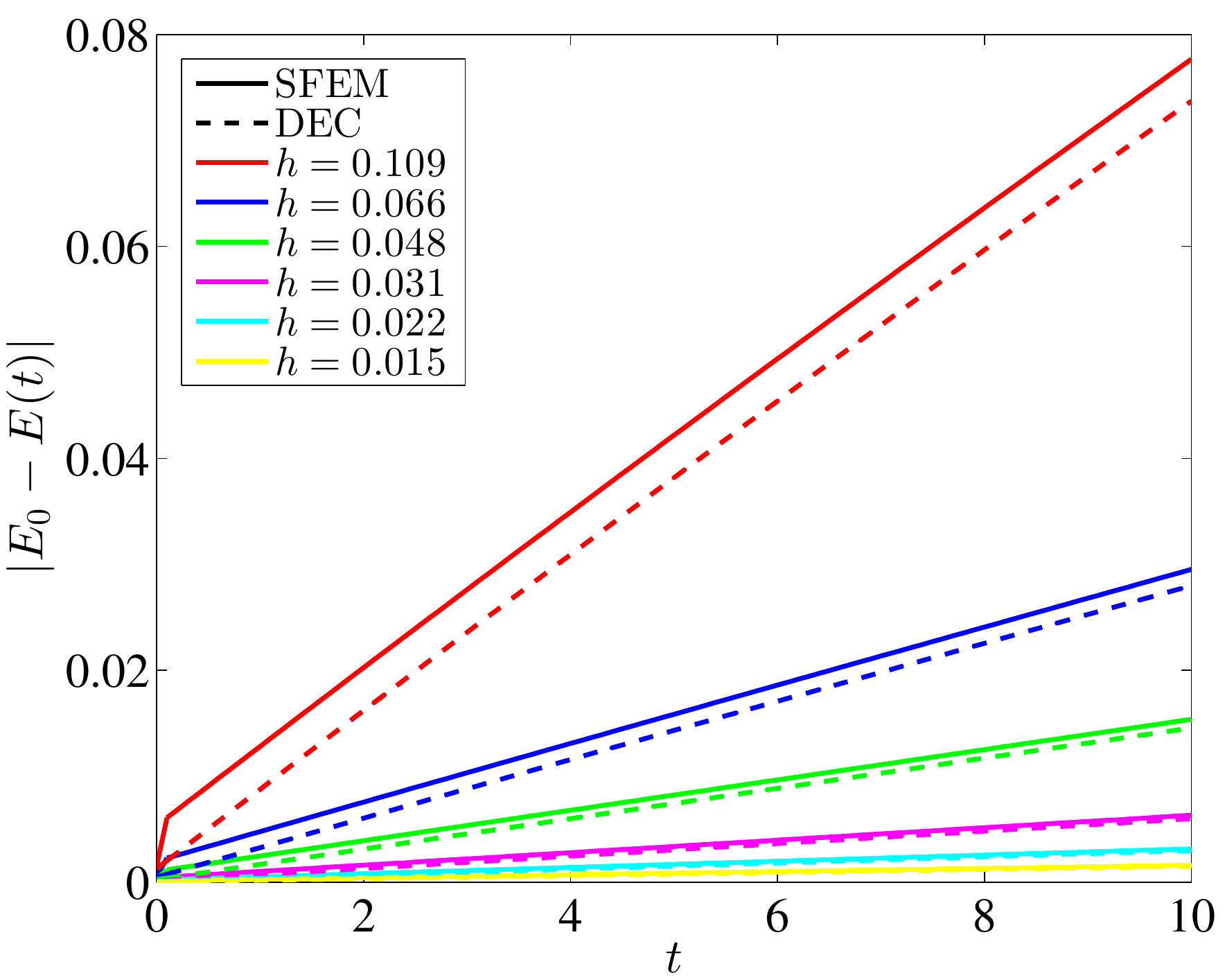} 
	\end{minipage}
	\hfill
	\begin{minipage}{0.49\textwidth}
		\includegraphics[width=\textwidth]{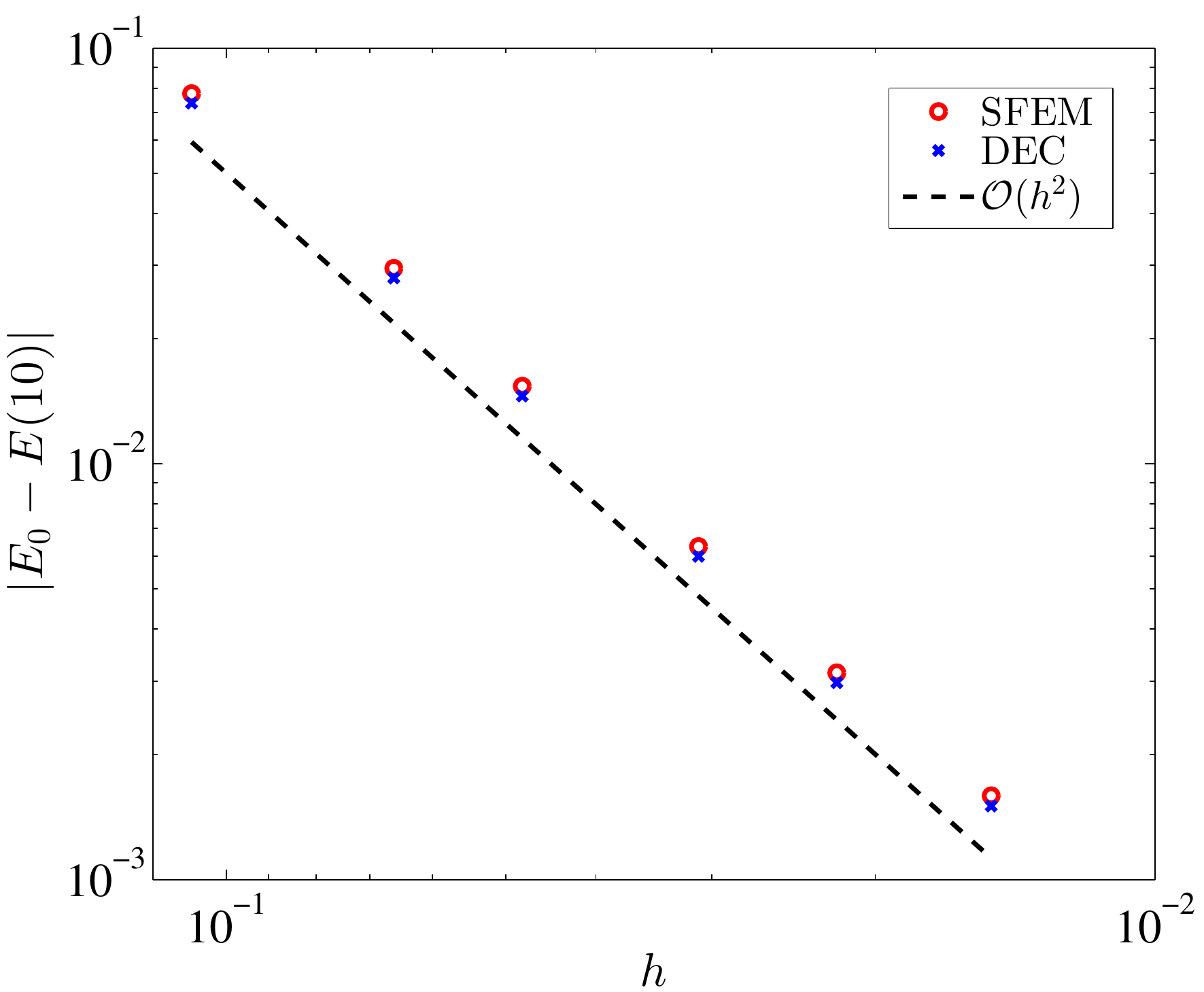}
	\end{minipage}
	\caption{Streamlines of stationary solution on a sphere together with the error in the kinetic energy and the experimental order of convergence (EOC) for different mesh sizes $h$ 
            (maximum circumcircle diameter of all triangles) at time $t = 10$ for both numerical methods. $E_0$ denotes the exact kinetic energy. The timestep is $\tau = 0.1$.}
	\label{fig2}
\end{figure}

We first consider a stationary solution on a sphere, with $\psi_0(\xb) = z$ and $\vb_0(\xb) = (y, -x,0)^T$ with coordinates $(x,y,z) \in \R^3$. Fig. \ref{fig2} $\;$ shows the streamlines for the rotating flow together with the computed errors for the kinetic energy. 
The results essentially show second order convergence for both methods, the DEC approach for the surface Navier-Stokes equation and the surface finite element method (SFEM) for the vorticity-stream function equation.
 
\subsection{Geometric interaction}

\begin{figure}[t]
	\begin{minipage}{0.1\textwidth}
		\includegraphics[width=\textwidth]{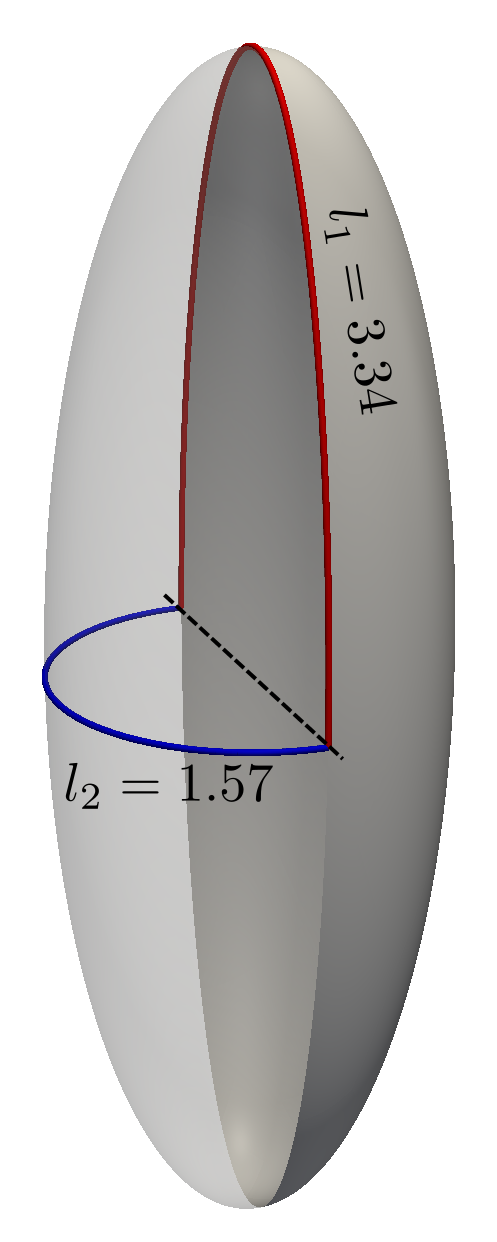}
	\end{minipage}
	\hfill
	\begin{minipage}{0.78\textwidth}
		\includegraphics[width=0.16\textwidth]{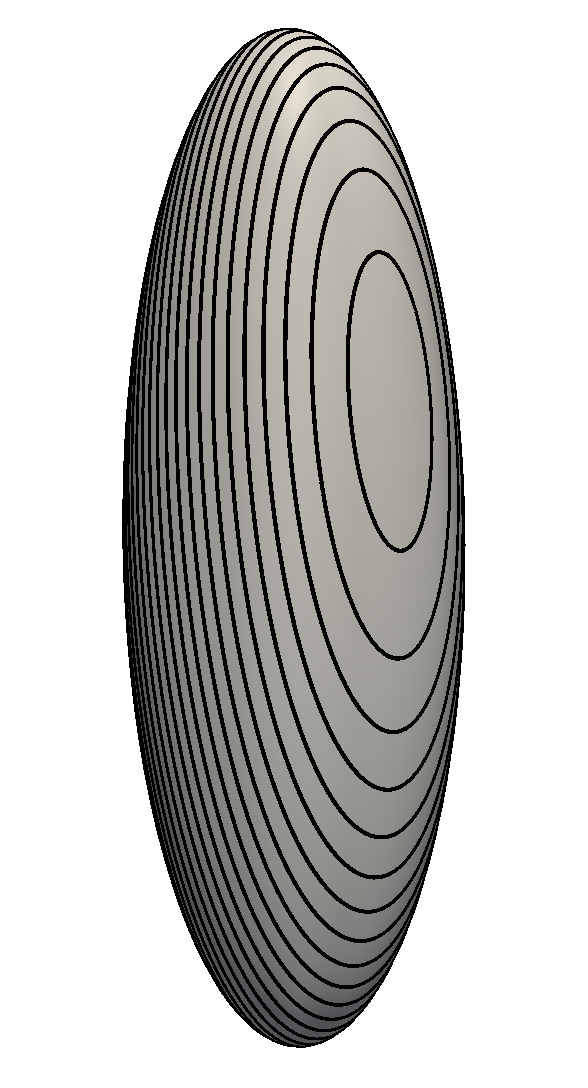}
		\includegraphics[width=0.16\textwidth]{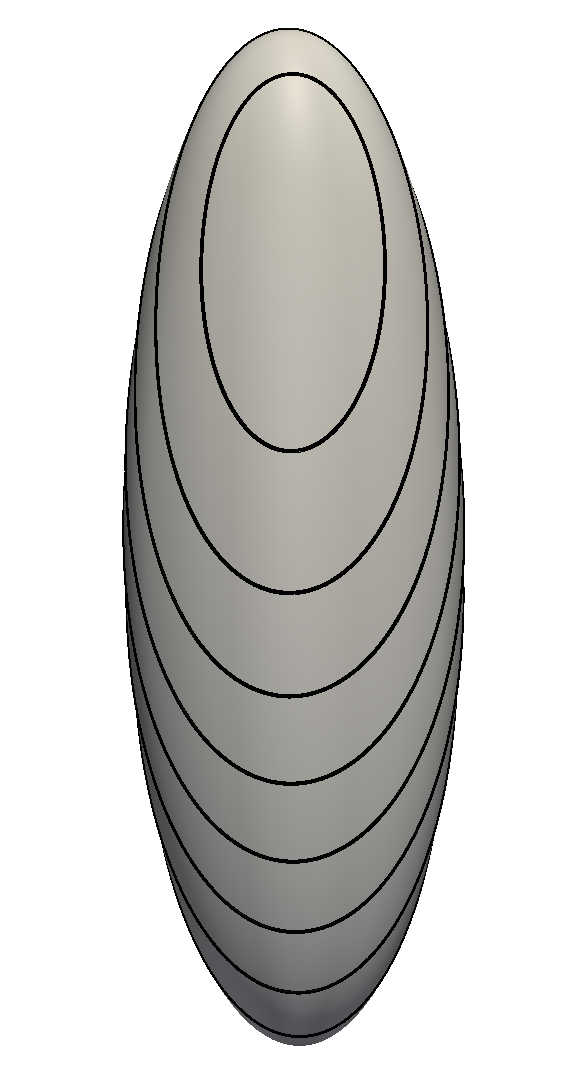}
		\includegraphics[width=0.16\textwidth]{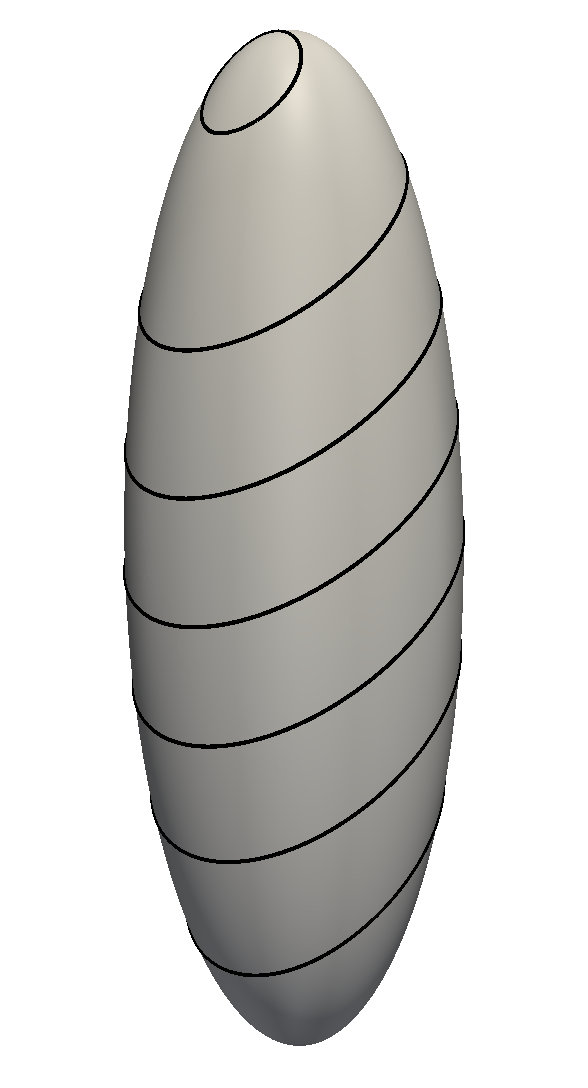}
		\includegraphics[width=0.16\textwidth]{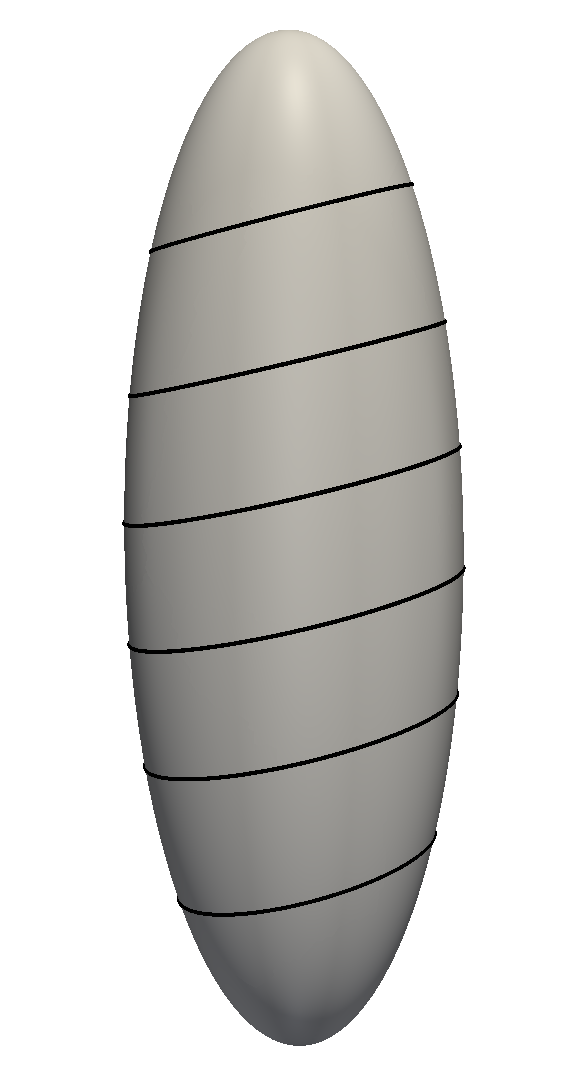}
		\includegraphics[width=0.16\textwidth]{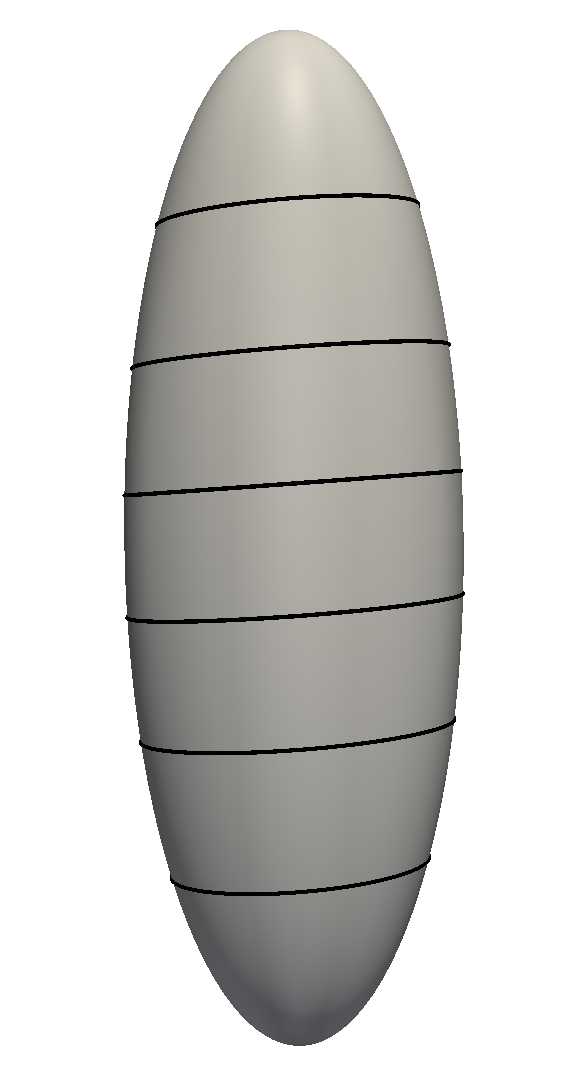}
		\includegraphics[width=0.16\textwidth]{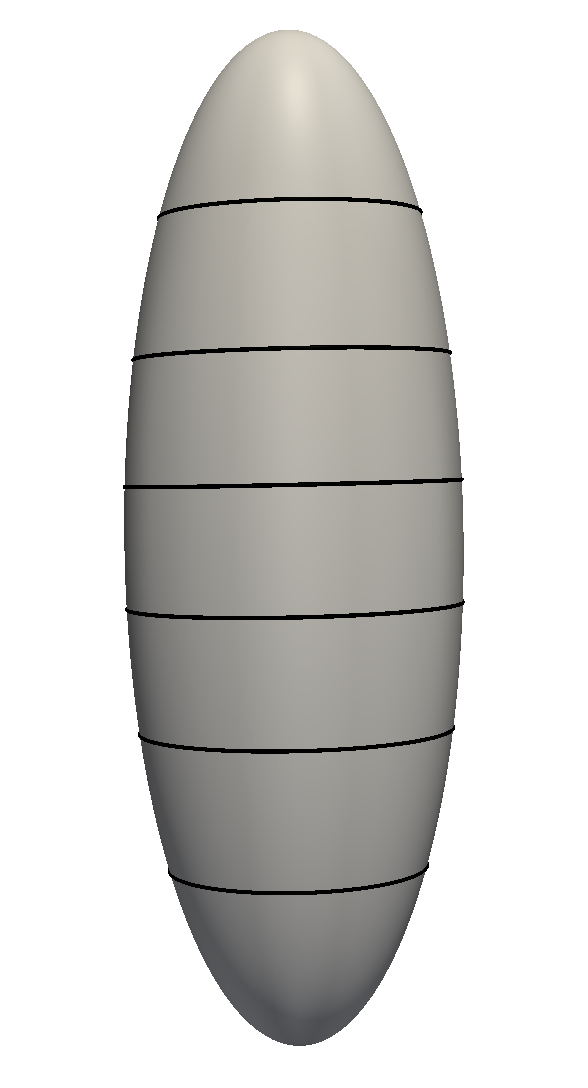}
	\end{minipage}\\
	\begin{minipage}{0.49\textwidth}
		\includegraphics[width=\textwidth]{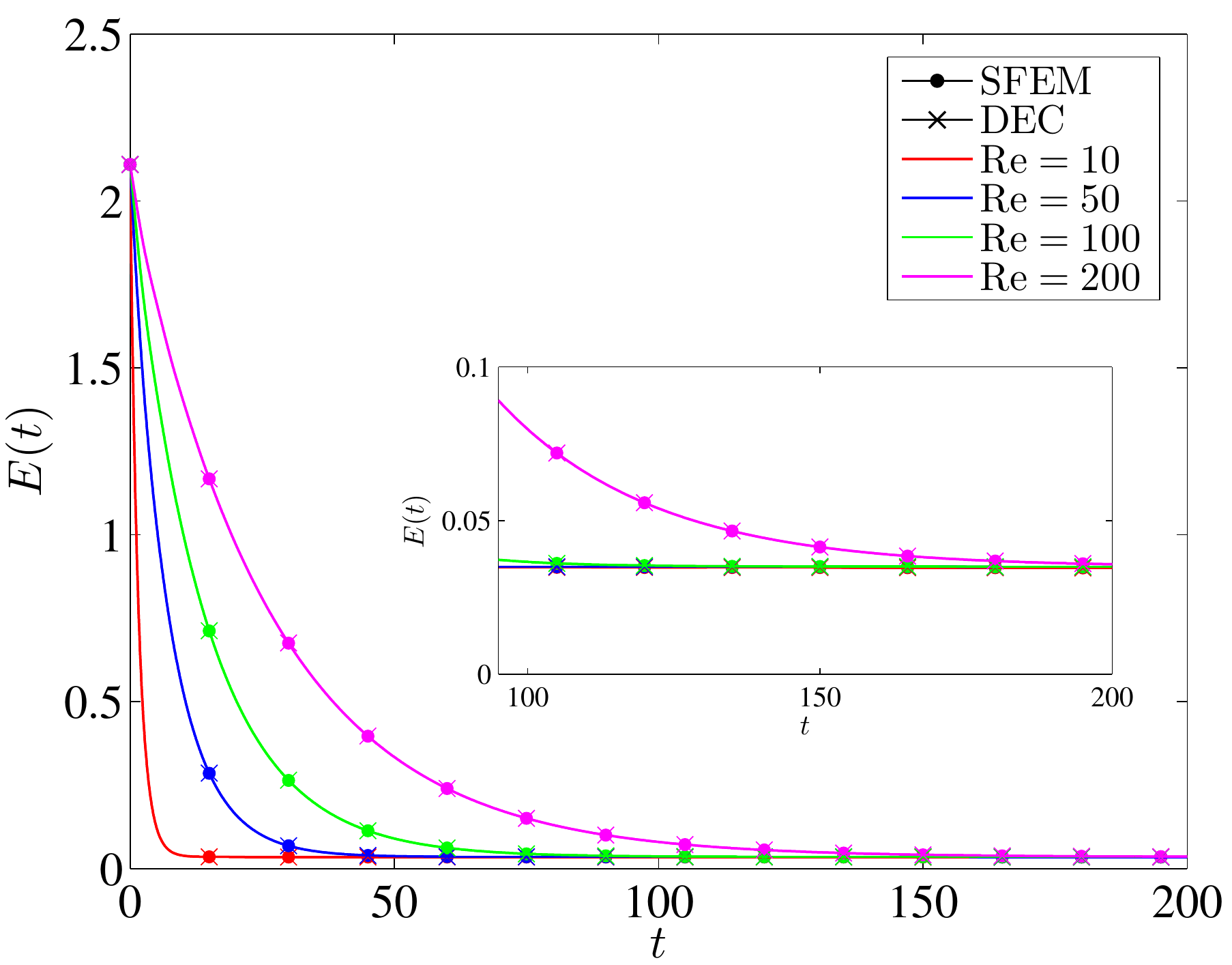}
	\end{minipage}
	\begin{minipage}{0.49\textwidth}
		\begin{tikzpicture}
			\node (vortexHeightPlot) at (0,0) {\includegraphics[width=\textwidth]{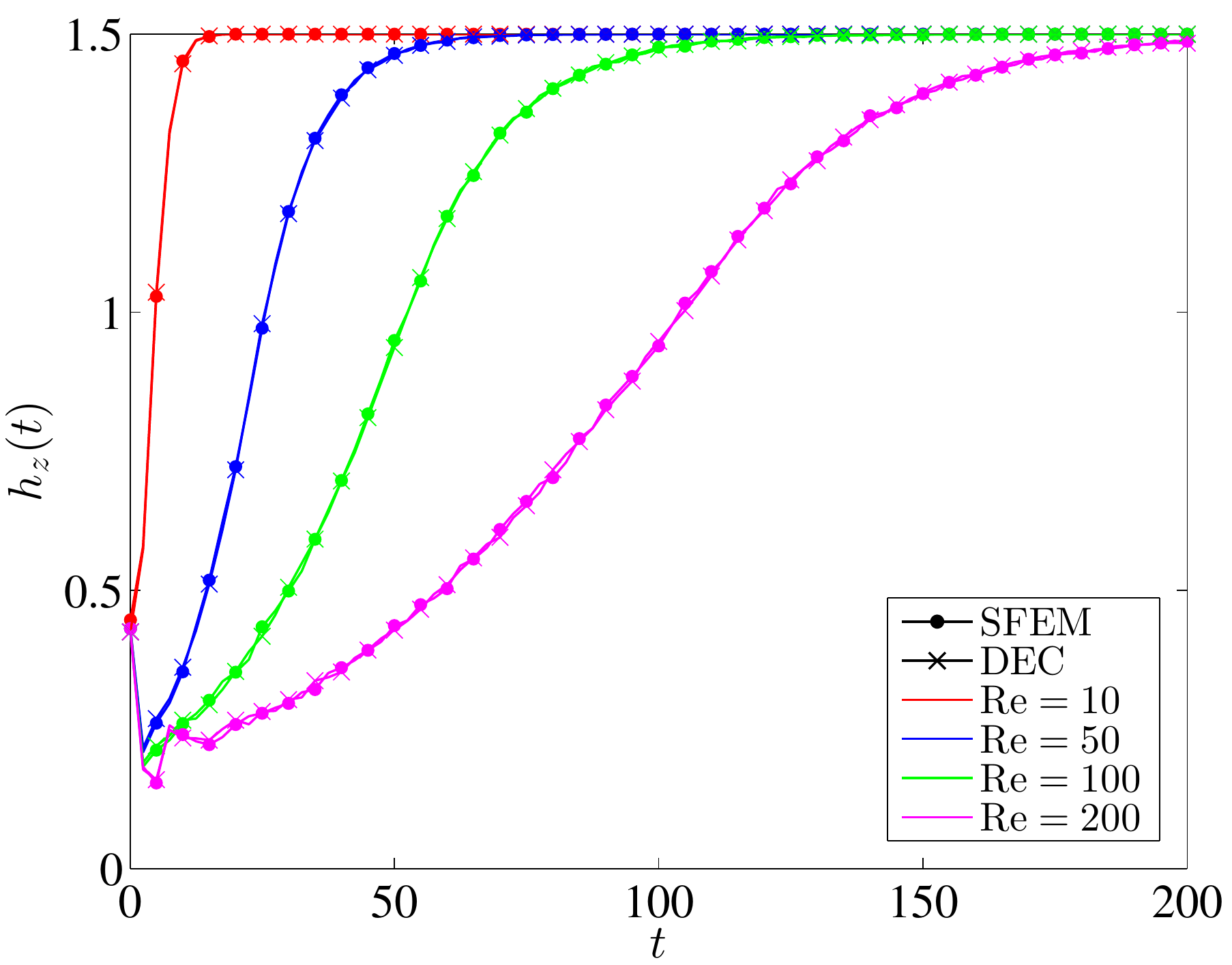}};
			\node (schematicHeight) at (1.9,0.5) {\includegraphics[width=0.18\textwidth]{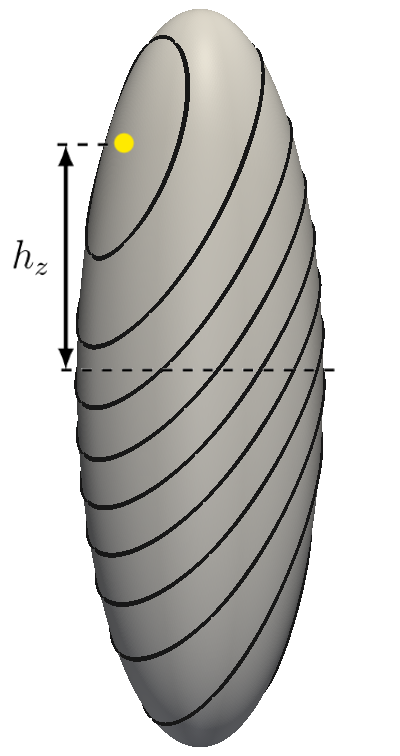}};
		\end{tikzpicture}
	\end{minipage}
	\caption{(top) Distances on the ellipsoid together with the streamlines at $t = 0, 4, 8, 12, 16$ and $20$ for the rotating flow. Results are shown for $\text{Re}= 10$. (bottom) Kinetic energy over time and the height for the upper vortex over time for both numerical approaches and various $\text{Re}$.}
	\label{fig3}
\end{figure}

As already analyzed in detail in \cite{Reutheretal_MMS_2015} the vortices in the flow, in the considered case two $+1$ defects, repel each other and are attracted by regions of high Gaussian curvature. We first consider an ellipsoid, represented by the level-set function $e(\xb) = (x/a)^2 + (y/b)^2 + (z/c)^2$, with $(x,y,z) \in \R^3$, $a = b = 0.5$ and $c = 1.5$. We consider the initial solutions $\psi_0(\xb) = y + 0.1 z$ and $\vb_0(\xb) = \RotSurf \psi_0(\xb)$ and use a timestep $\tau = 0.1$. Fig. \ref{fig3} shows the geometric properties, the streamlines at various times for $\text{Re} = 10$ as well as the kinetic energy over time and the position of one vortex over time for both methods and various $\text{Re}$. 
The flow converges to a stationary solution with the vortices located at the high Gaussian curvature regions. However, these positions also favors the long range interaction between the vortices as they maximize their distance. We thus cannot argue on a geometric interaction. The time to reach the stationary solution strongly depends on $\text{Re}$, the lower $\text{Re}$ the faster it is reached.  
 
\begin{figure}[t]
	\centering
	\begin{minipage}{0.6\textwidth}
		\begin{minipage}{0.46\textwidth}
			\includegraphics[width=\textwidth]{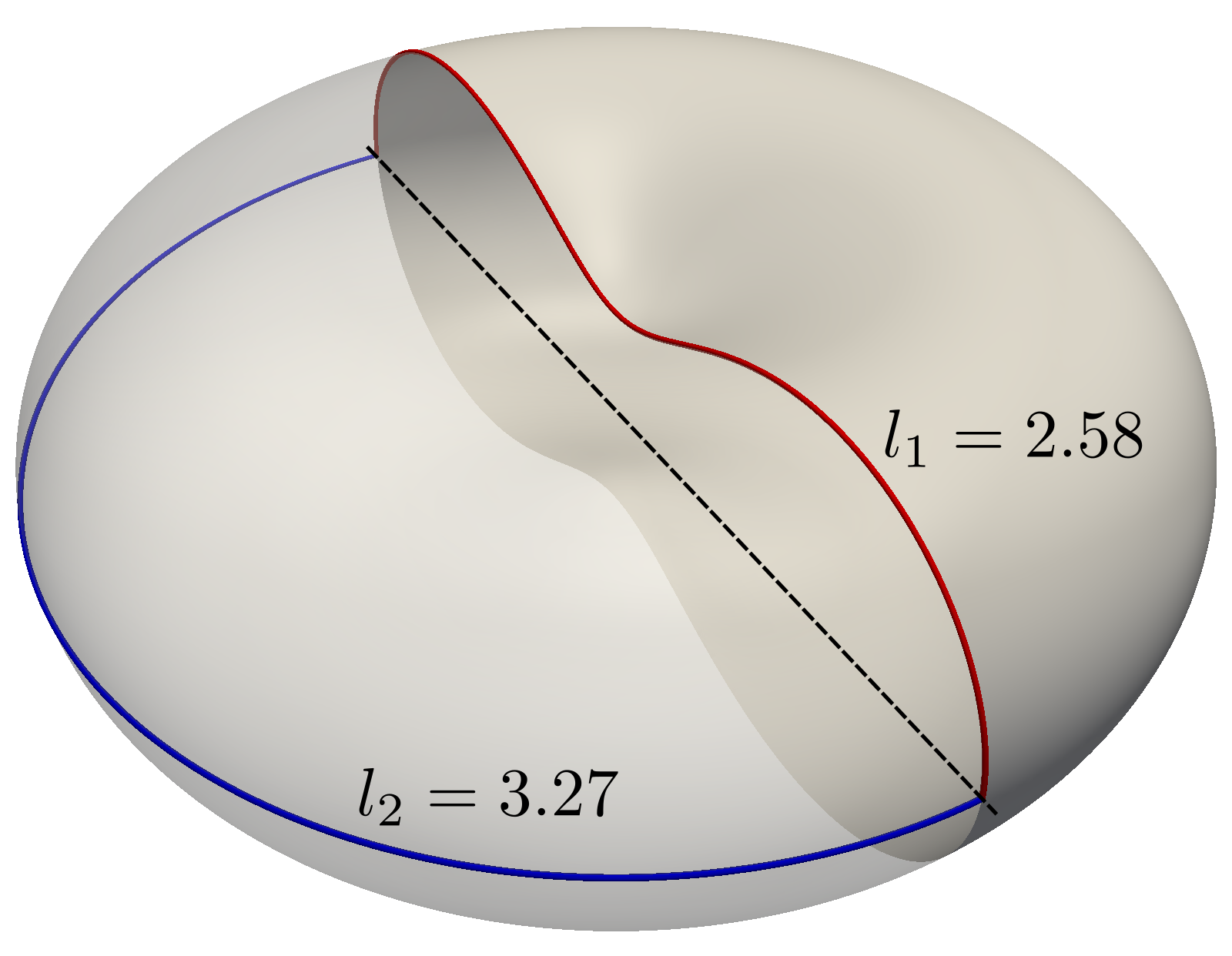}
		\end{minipage}
		\begin{minipage}{0.46\textwidth}
			\includegraphics[width=\textwidth]{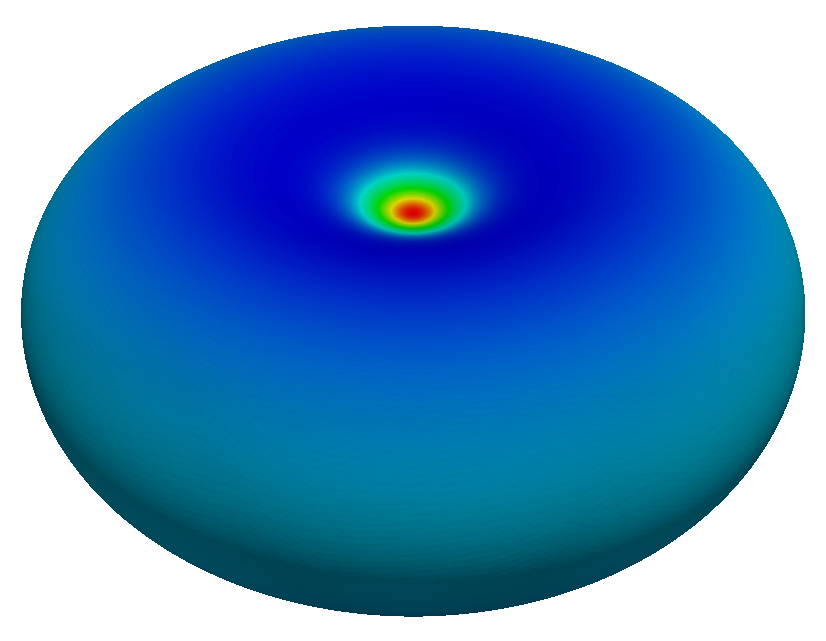}
		\end{minipage}
		\begin{minipage}{0.06\textwidth}
			\includegraphics[width=\textwidth]{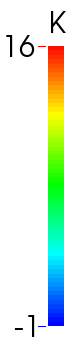}
		\end{minipage}
	\end{minipage}\\
	\begin{minipage}{\textwidth}
		\includegraphics[width=0.19\textwidth]{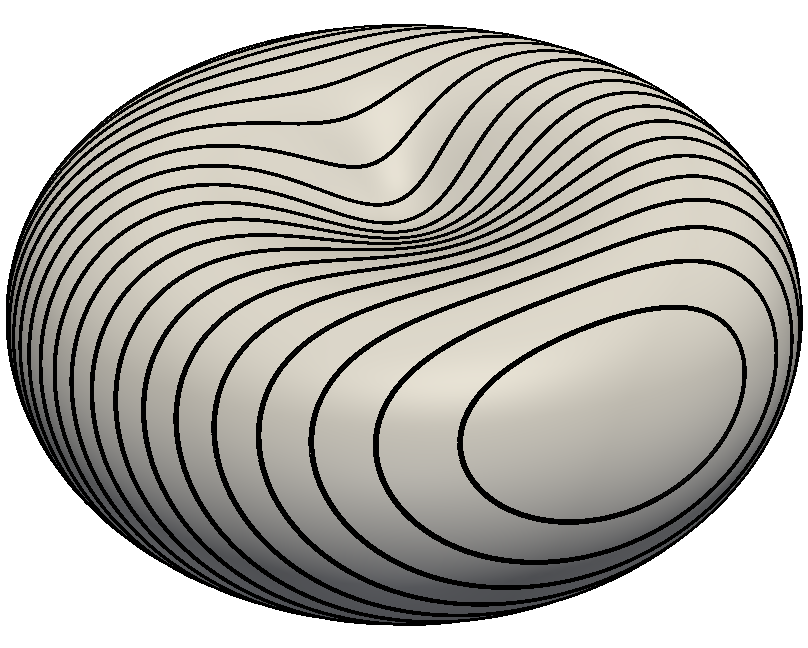}
		\includegraphics[width=0.19\textwidth]{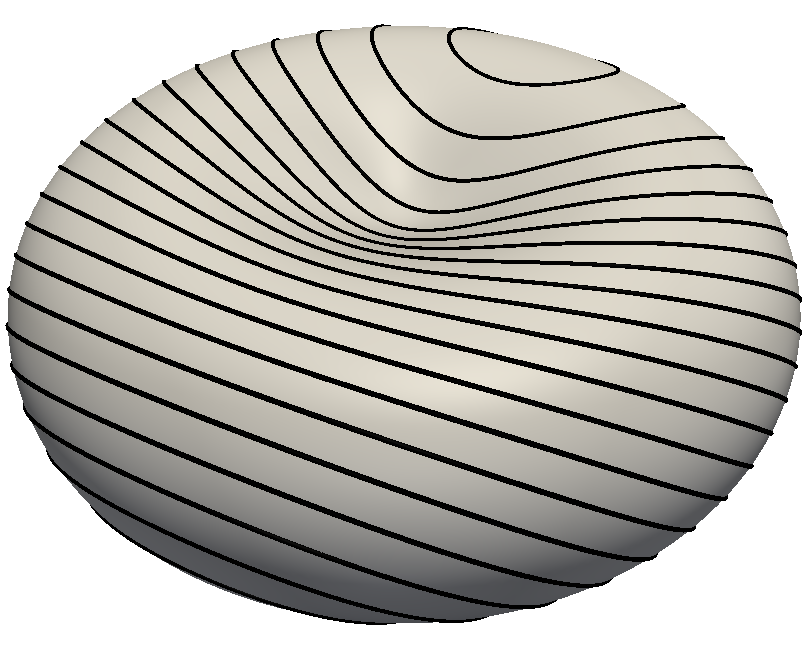}
		\includegraphics[width=0.19\textwidth]{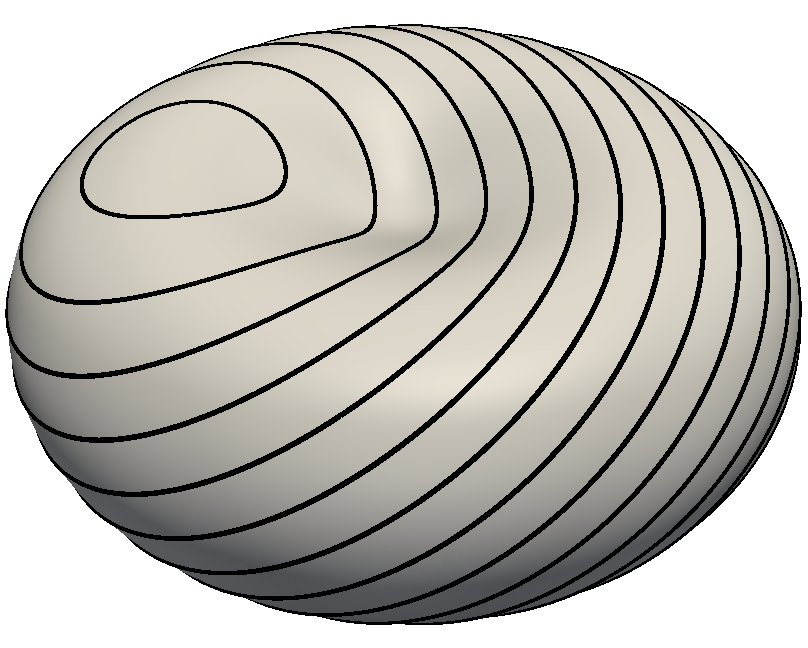}
		\includegraphics[width=0.19\textwidth]{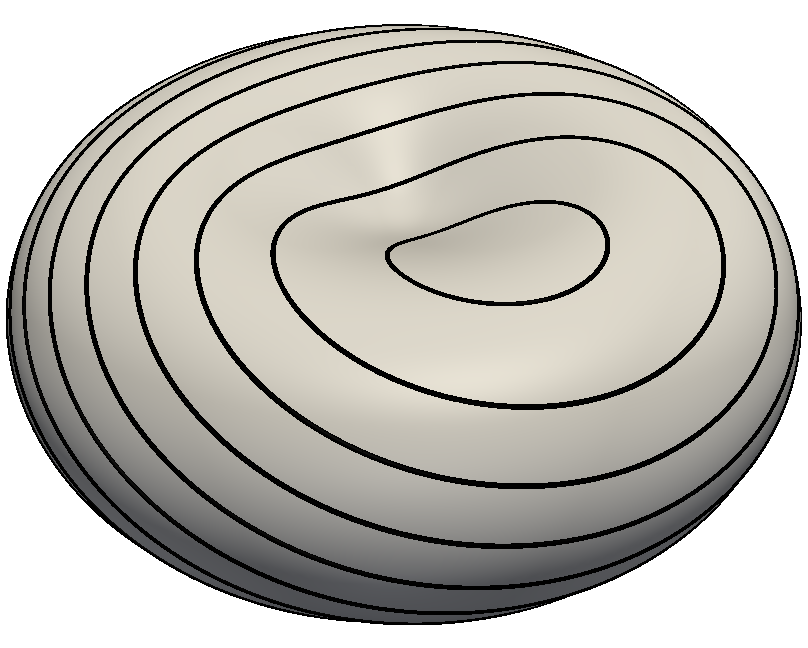}
		\includegraphics[width=0.19\textwidth]{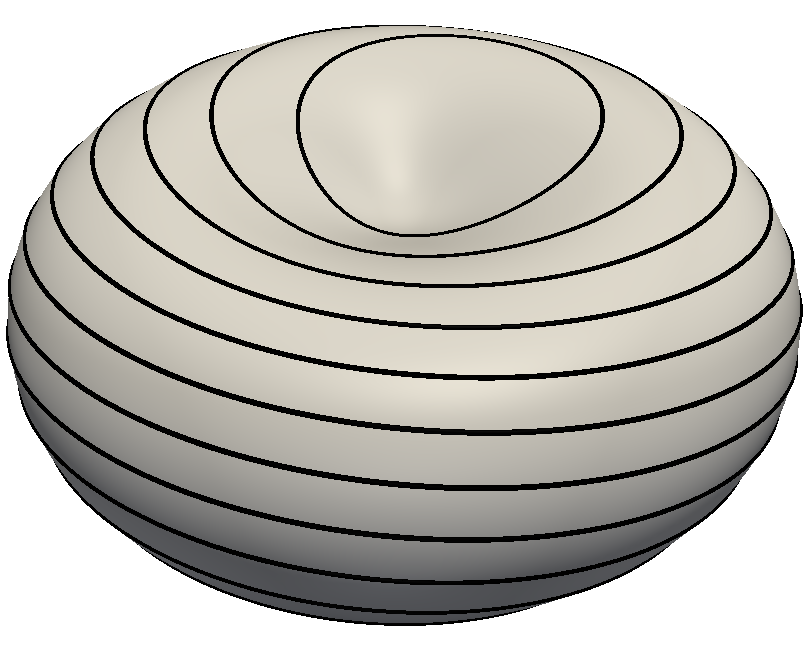}\\
		\includegraphics[width=0.19\textwidth]{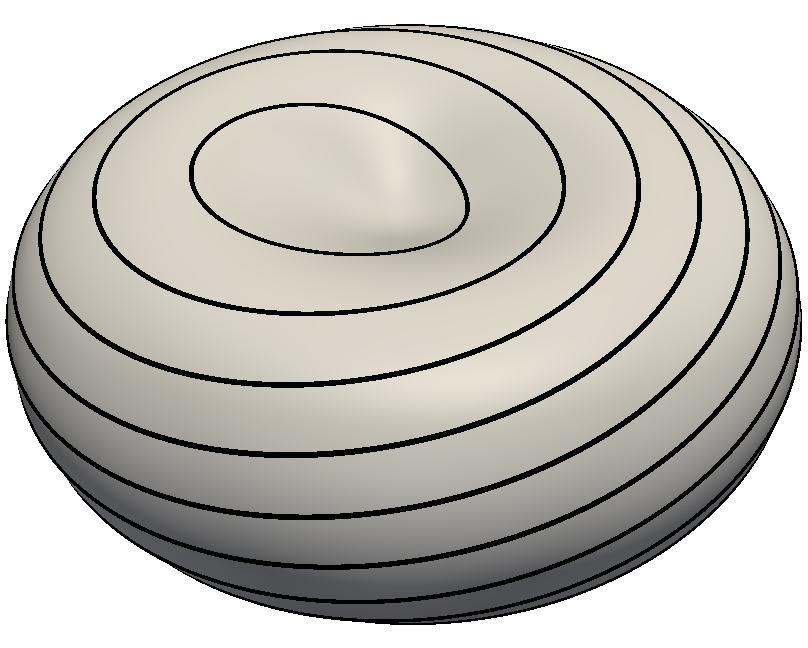}
		\includegraphics[width=0.19\textwidth]{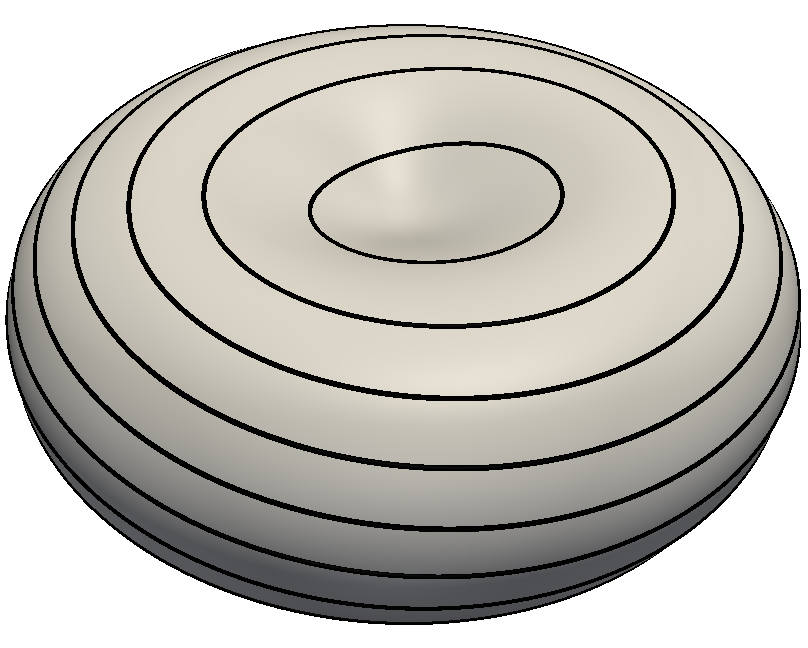}
		\includegraphics[width=0.19\textwidth]{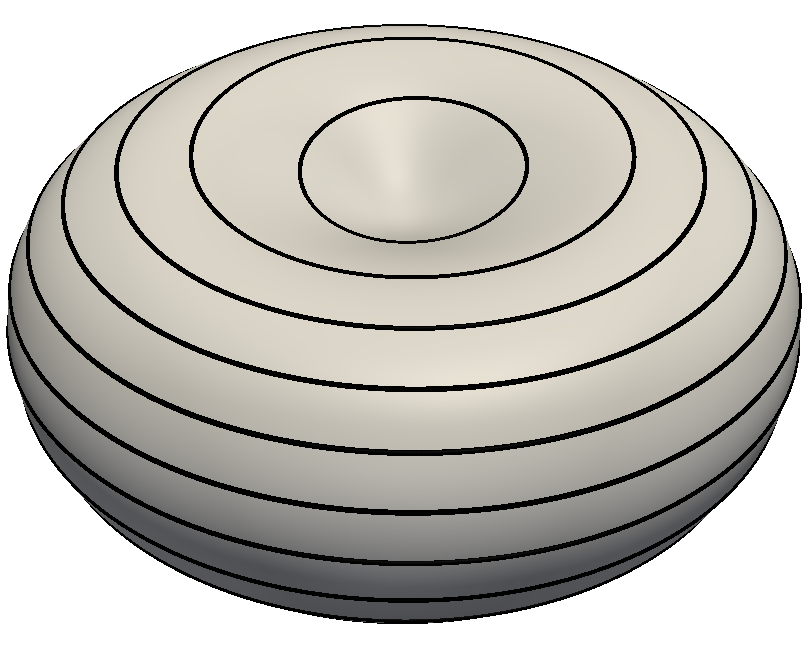}
		\includegraphics[width=0.19\textwidth]{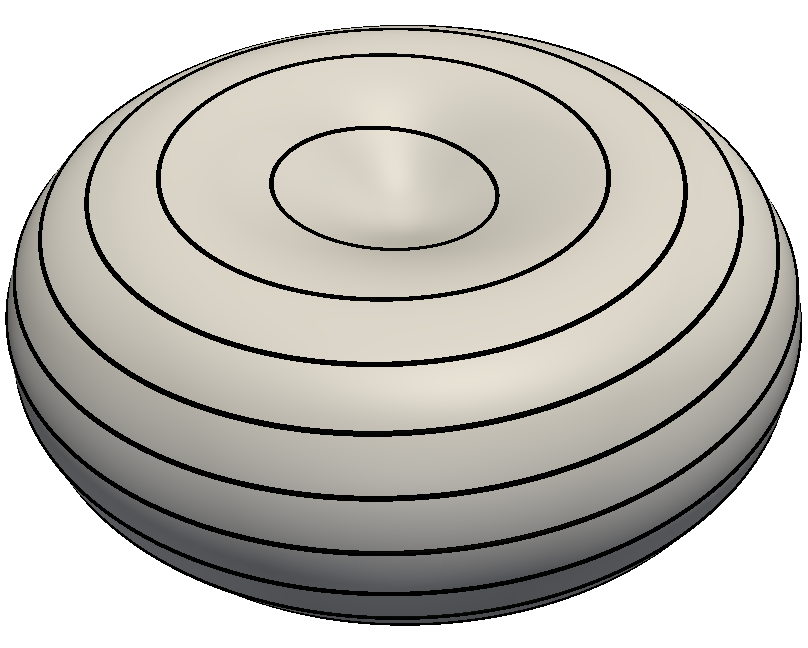}
		\includegraphics[width=0.19\textwidth]{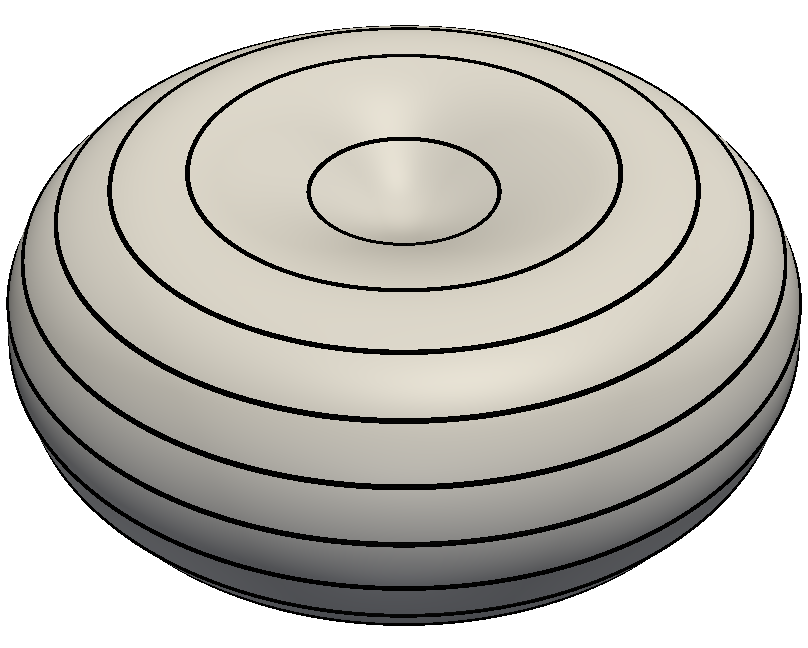}
	\end{minipage} \\
	\begin{minipage}{\textwidth}
		\begin{minipage}{0.45\textwidth}
			\includegraphics[width=\textwidth]{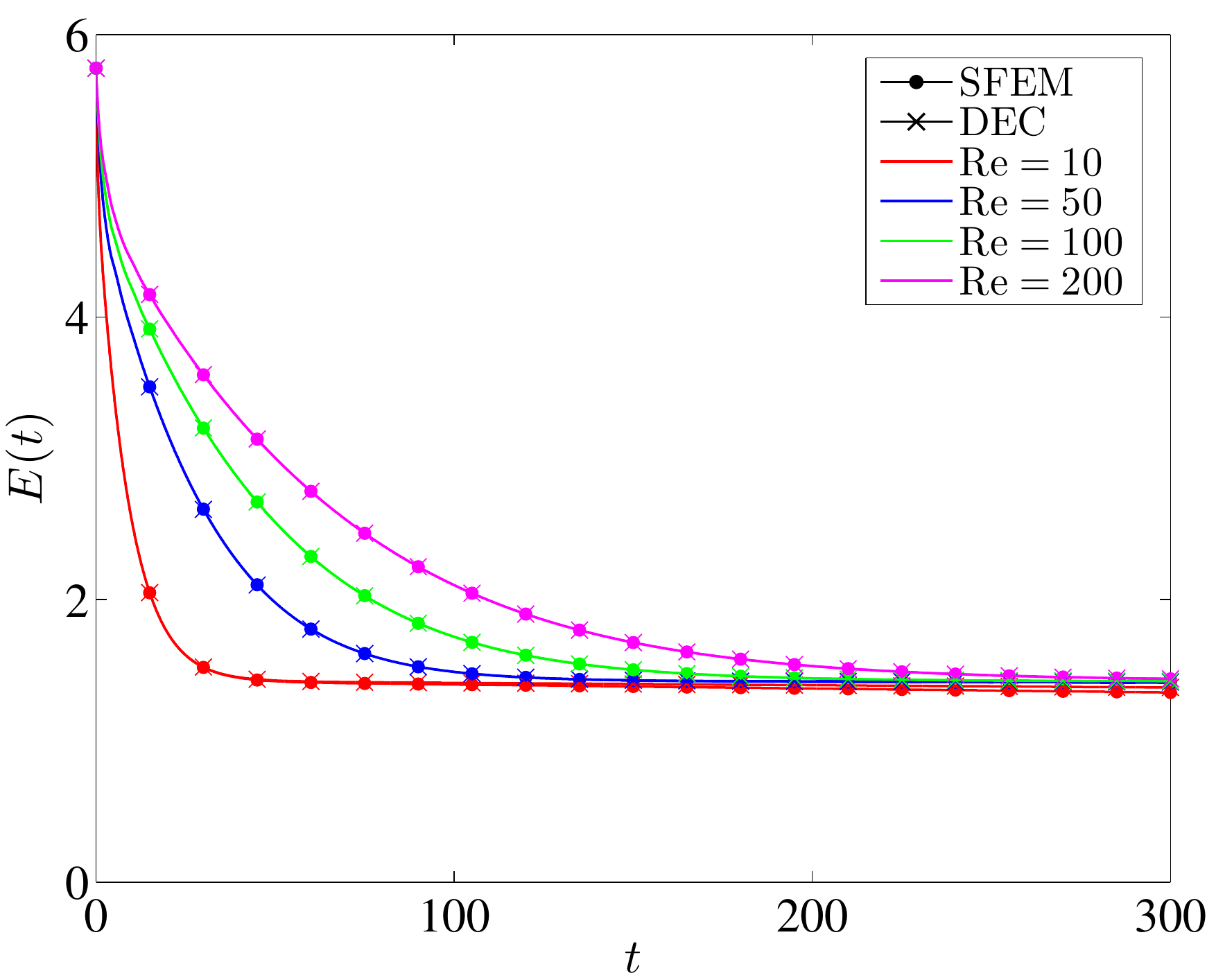}
		\end{minipage}
		\begin{minipage}{0.26\textwidth}
			\includegraphics[width=\textwidth]{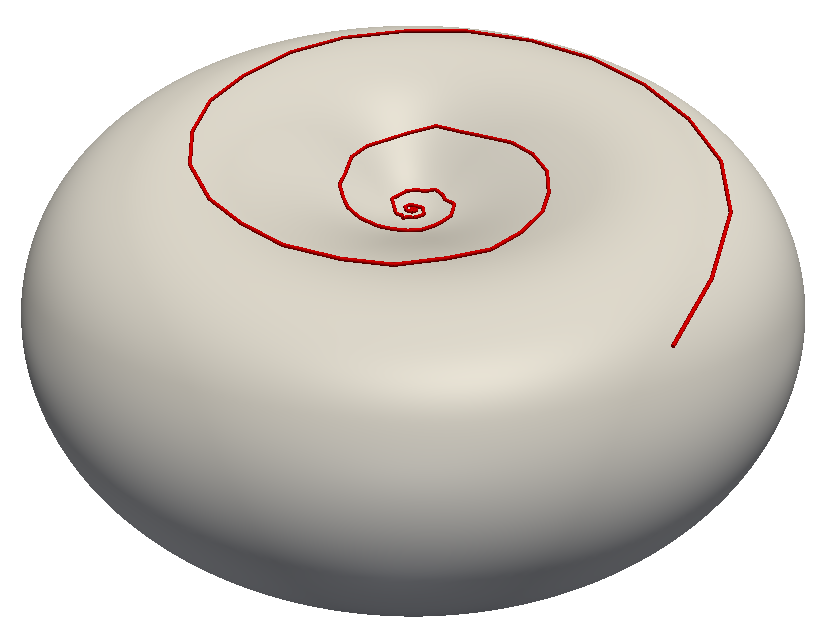}
		\end{minipage}
		\begin{minipage}{0.26\textwidth}
			\includegraphics[width=\textwidth]{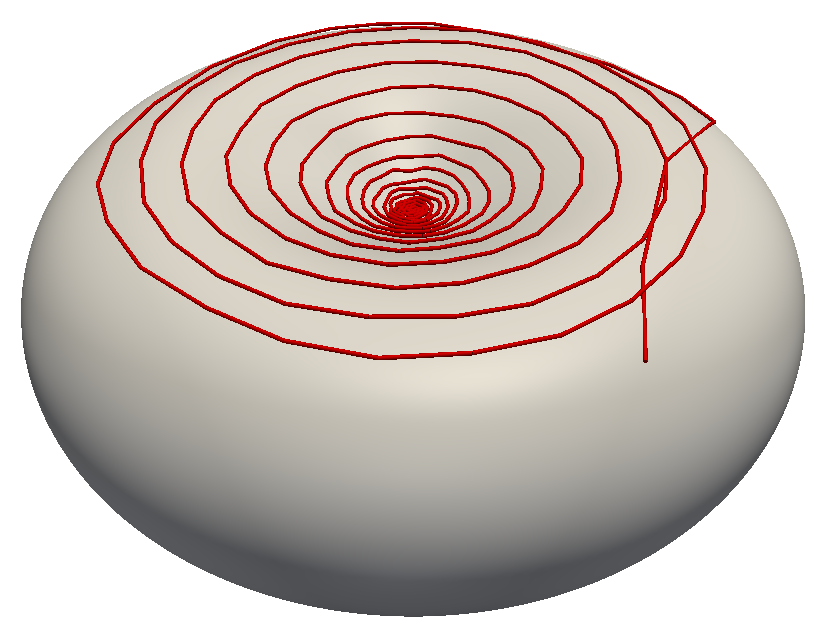} 
		\end{minipage}
	\end{minipage}
	\caption{(top) Distances on the biconcave shape together with the Gaussian curvature. (middle) Streamlines at $t = 0, 7, 14, 21, 28, 35, 42, 49, 56$ and $200$ (left to right, top to bottom) for the rotating flow. Results are shown for $\text{Re}= 10$. (bottom) Kinetic energy over time for both numerical approaches and various $\text{Re}$ together with two examples for the vortex trajectories for $\text{Re} = 10$ (left) and $\text{Re} = 100$ (right).}
	\label{fig4}
\end{figure}

The second example considers a biconcave shape, represented by the level-set function $e(\xb) = (a^2 + x^2 + y^2 + z^2)^3 - 4 a ^2 ( y^2 + z^2) - c^4$, with $(x,y,z) \in \R^3$, $a = 0.72$ and $c = 0.75$. We consider the initial solutions $\psi_0(\xb) = y +z$ and $\vb_0(\xb) = \RotSurf \psi_0(\xb)$ and use a timestep $\tau = 0.1$. Fig. \ref{fig4} shows the geometric properties together with the trajectories of one vortex for different $\text{Re}$, the streamlines at various times, a plot of the Gaussian curvature and the kinetic energy over time. 
Again the flow converges to a stationary solution with the vortices located at the high Gaussian curvature regions. Here the location of the vortices clearly is a result of the geometric interaction, as their distance is not maximized. Again the time to reach the stationary solution strongly depends on $\text{Re}$, the lower $\text{Re}$ the faster it is reached.  

\subsection{Surfaces with genus $g(\surf) \neq 0$}
As mentioned above we will consider nontrivial solutions with $\DivSurf \vb =0$ and $\RotSurf \vb = 0$. The vorticity-stream function formulation in eqs. \eqref{eq:vor1} and \eqref{eq:vor2} is based on the Hodge decomposition of the velocity field $\vb$ which can be written as 
\begin{equation}
	\vb = \vb^{div} + \vb^{rot} + \vb^{harm}
	\label{eq:hodge_decomposition}
\end{equation}
on a general surface $\surf$ with a divergence free vector field $\vb^{div}$, a curl free vector field $\vb^{rot}$ and a divergence as well as curl free vector field $\vb^{harm}$. The first two parts are usually rewritten as $\vb^{div} = \RotSurf\psi$ and $\vb^{rot} = \GradSurf\Phi$ with scalar functions $\psi$ and $\Phi$. Since we require incompressibility of $\vb$ one can easily verify that the curl free part $\vb^{rot}$ vanishes identically. Furthermore, on spherical surfaces ($g(\surf) = 0$) we can drop the harmonic part since it is not possible to write a vector field that is divergence and curl free except of the zero vector field. Finally, this leads to the substitution $\vb = \RotSurf\psi$ which is used in the vorticity-stream function approach in the prior sections. On surfaces with $g(\surf) \neq 0$ the situation changes and the harmonic part $\vb^{harm}$ does not vanish generally. To demonstrate this property we use the torus which has genus $g(\surf) =1$. A torus can be described by the levelset function $e(\xb) = ( \sqrt{x^2+z^2} - R )^2 + y^2 - r^2$, with $(x,y,z) \in \R^3$, major radius $R$ and minor radius $r$. Throughout this section we use $R = 2$ and $r=0.5$. Let $\phi$ and $\theta$ denote the standard parametrization angles  on the torus. Then, the two basis vectors can be written as $\partial_\phi\xb$ as well as $\partial_\theta\xb$ and read in Cartesian coordinates $\partial_\phi\xb = (-z, 0, x)$ as well as $\partial_\theta\xb = (-\frac{xy}{\sqrt{x^2+z^2}},  \sqrt{x^2+z^2} - 2, -\frac{yz}{\sqrt{x^2+z^2}})$ 
%\begin{equation*}
%	\partial_\phi\xb = \left( \begin{matrix} -z \\ 0 \\ x \end{matrix} \right) \qquad \mbox{as well as} \qquad \partial_\theta\xb = \left( \begin{matrix} -\frac{xy}{\sqrt{x^2+z^2}} \\ \sqrt{x^2+z^2} - 2 \\ -\frac{yz}{\sqrt{x^2+z^2}} \end{matrix} \right) \formComma
%\end{equation*}
which are schematically shown in \autoref{fig:torus:1}. We find two (linear independent) harmonic vector fields on the torus
\begin{align*}
	\vb_\phi^{harm} &= (4+\cos(\theta))^{-2} \partial_\phi\xb = \frac{1}{4\left(x^2+z^2\right)} \partial_\phi\xb \\
	\vb_\theta^{harm} &= (4+\cos(\theta))^{-1} \partial_\theta\xb = \frac{1}{2\sqrt{x^2+z^2}} \partial_\theta\xb
\end{align*}
written in local and Cartesian coordinates, respectively, and shown in \autoref{fig:torus:2}. One can easily verify that $\DivSurf\vb_\phi^{harm} = \RotSurf\vb_\phi^{harm} = 0$ as well as $\DivSurf\vb_\theta^{harm} = \RotSurf\vb_\theta^{harm} = 0$.

To start with, consider the vector field $\vb = \partial_\phi\xb$, which has zero divergence and non-zero curl. The Hodge decomposition eq. \eqref{eq:hodge_decomposition} leads to $\vb^{rot}=\vb^{harm}=0$. In that case the substitution $\vb = \RotSurf(\psi)$ holds. The stream function $\psi$ of the vector field $\partial_\phi\xb$ can then be analytically written in local coordinates as $\psi = -\frac{1}{4}\sin(\theta) + \theta - \pi$. The linear contribution causes a discontinuity at $\theta = 2\pi$, which is shown in \autoref{fig:torus:1} together with the streamlines of $\partial_\phi\xb$ (contour lines of $\psi$). Solving the surface Navier-Stokes equation eqs. \eqref{eq1} and \eqref{eq2} directly circumvents the discontinuities.

In the next example we use the mean of the two harmonic vector fields as initial condition, i.e.
$\vb_0(\xb) = \frac{1}{2}( \vb_\phi^{harm} + \vb_\theta^{harm} )$. 
By considering the vorticity-stream function approach we have the initial conditions $\phi_0 = \psi_0 = 0$ and thus only the trivial solution. However, solving the surface Navier-Stokes equation directly covers also the harmonic parts. \autoref{fig:torus:2} shows the numerical solution of $\vb$ with the DEC algorithm in which we used the timestep $\tau=0.1$ and $\text{Re}=10$. In this case the reached steady state is again a Killing vector field and is proportional to the basis vector $\partial_\phi\xb$. Interestingly, the curl of the vector field $\partial_\phi\xb$ does not vanish.

Other linear combinations of the two harmonic vector fields $\vb_\phi^{harm}$ and $\vb_\theta^{harm}$ as initial condition leads to the same steady state solution (up to a proportionality constant) except of $\vb_0(\xb) = \vb_\theta^{harm}$. In that case the vector field does not change its direction by symmetry and dissipates to zero. The results are shown in the energy plot in \autoref{fig:torus:2} which clearly shows the vanishing energy over time.

\begin{figure}[t]
	\centering
	\begin{minipage}{0.3\textwidth}
		\begin{tikzpicture}
			% torus shape
			\node (vortexHeightPlot) at (0,0) {\includegraphics[width=\textwidth]{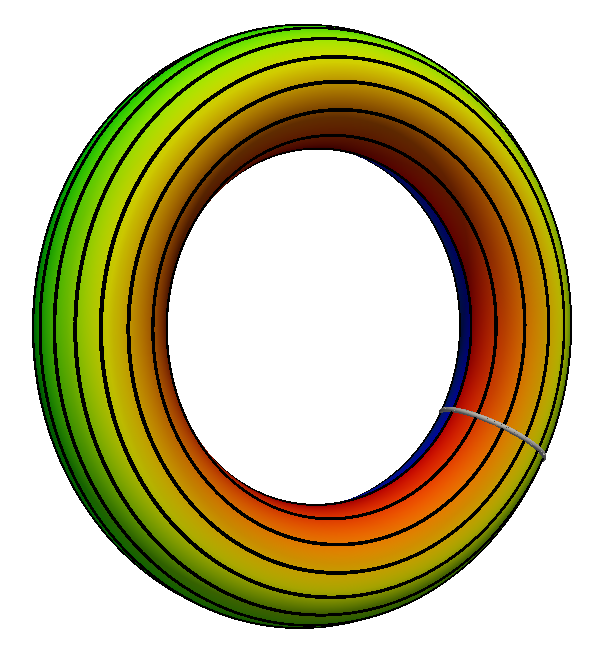}};
			% local coordinate u
			\draw[color=white,thick] (-1.34,0) to (-1.14,0) node {};
			\draw[color=white,arrows={-latex},thick] (-1.24,0) to [out=90,in=-110] (-1.15,0.5) node {};
			\draw (-1.1,-0.1) node[color=white,anchor=north] {\large$\partial_\phi\xb$};
			% local coordinate v
			\draw[color=white,thick] (0.25,-1.5) to (0.45,-1.48) node {};
			\draw[color=white,arrows={-latex},thick] (0.35,-1.49) to [out=84,in=-62] (0.25,-1.05) node {};
			\draw (0.3,-1.39) node[color=white,anchor=east] {\large$\partial_\theta\xb$};
		\end{tikzpicture}
	\end{minipage}
	\begin{minipage}{0.07\textwidth}
		\centering
		$\psi$\\
		\includegraphics[width=\textwidth]{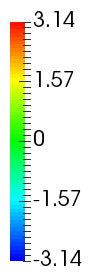}
	\end{minipage}
	\hspace{0.05\textwidth}
	\begin{minipage}{0.4\textwidth}
		\includegraphics[width=\textwidth]{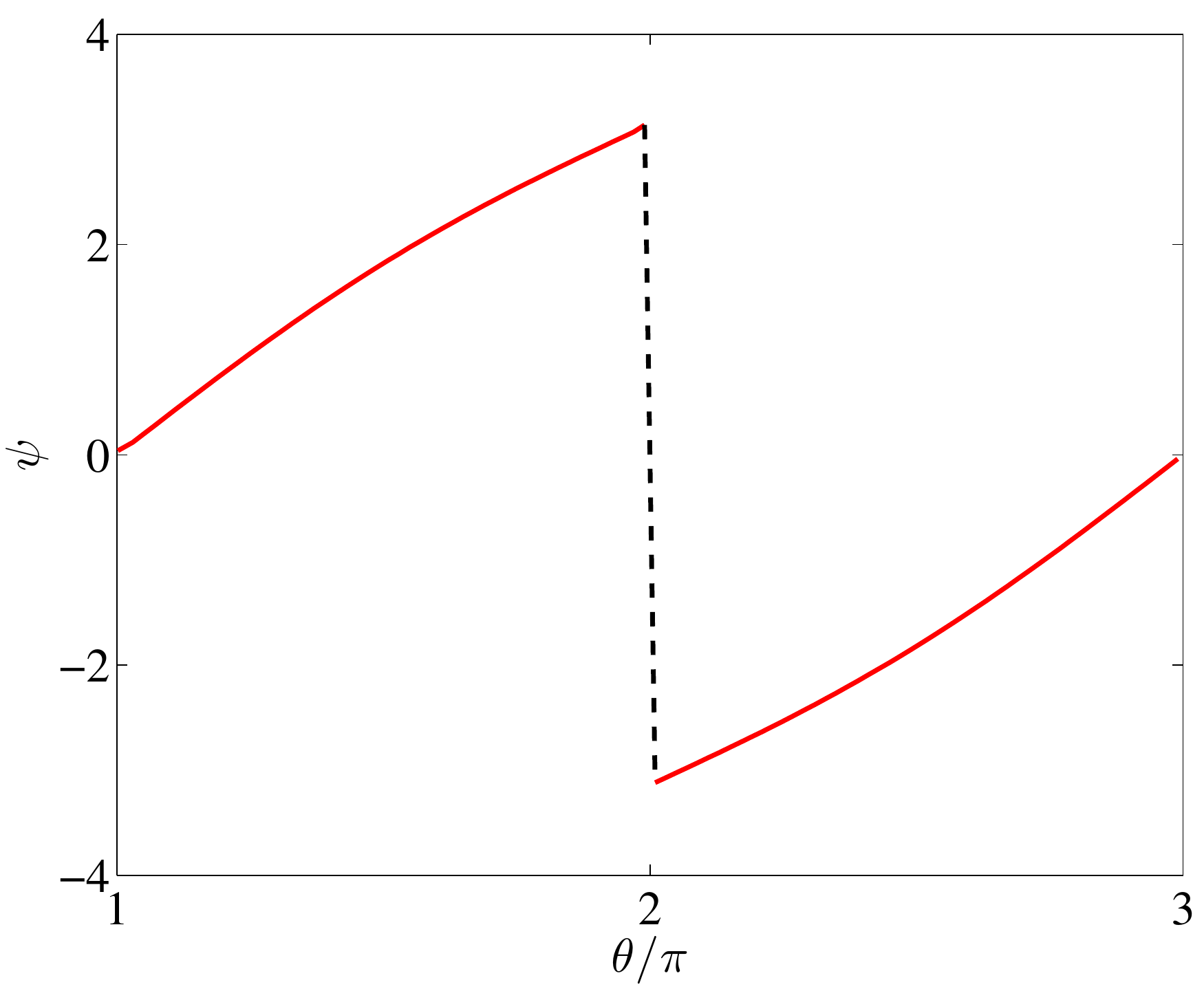}
	\end{minipage}
	\caption{(left) Streamlines and values of the discontinuous stream function $\psi$ to represent the velocity field $\vb = \partial_\phi\xb$ on the torus and the two basis vectors $\partial_\phi\xb$ and $\partial_\theta\xb$. (right) Plot of the stream function values over the gray contour line in the left figure.}
	\label{fig:torus:1}
\end{figure}
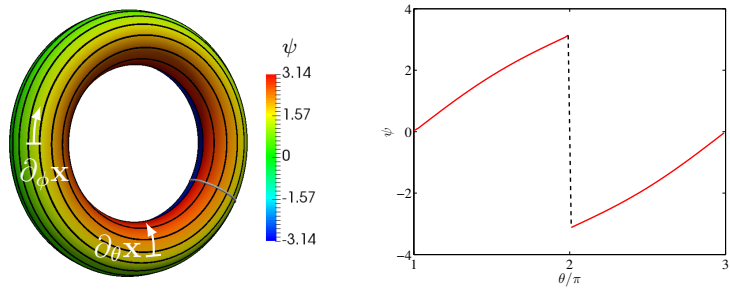
\begin{figure}[t]
	\centering
	\begin{minipage}{0.2\textwidth}
		\centering
		\includegraphics[width=\textwidth]{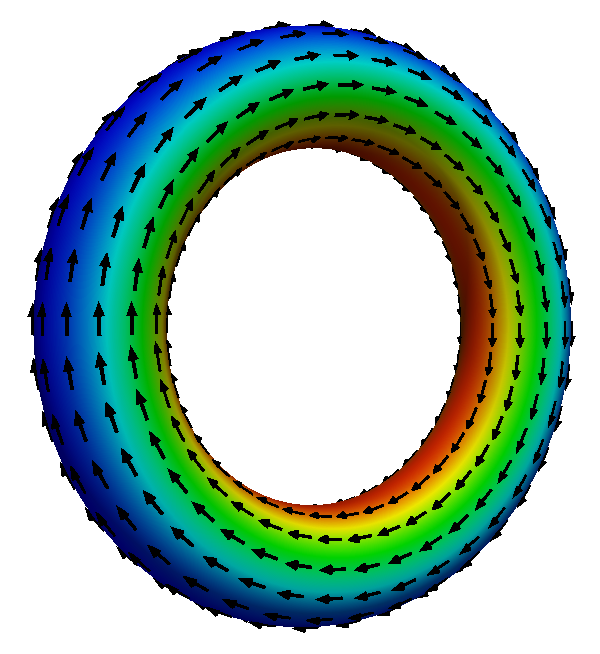}\\
		$\vb_\phi^{harm}$
	\end{minipage}
	\begin{minipage}{0.2\textwidth}
		\centering
		\includegraphics[width=\textwidth]{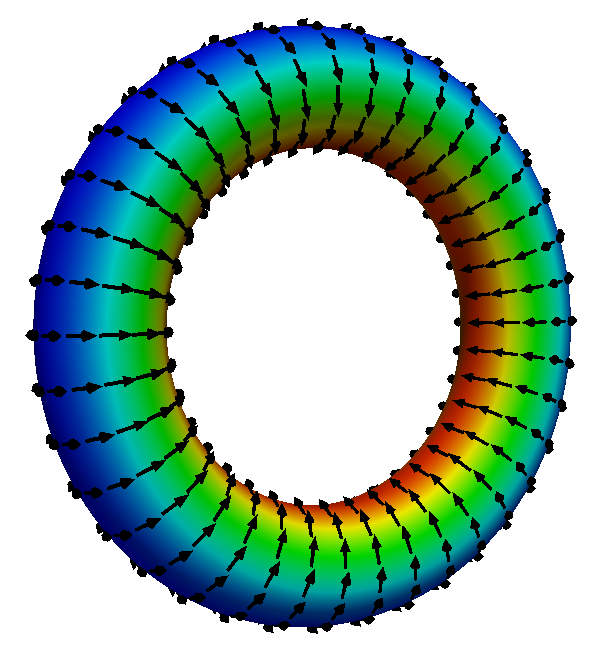}\\
		$\vb_\theta^{harm}$
	\end{minipage}
	\begin{minipage}{0.07\textwidth}
		\centering
		\includegraphics[width=\textwidth]{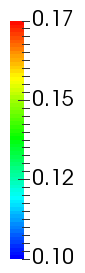}
	\end{minipage}
	\hspace{0.05\textwidth}
	\begin{minipage}{0.4\textwidth}
		\includegraphics[width=\textwidth]{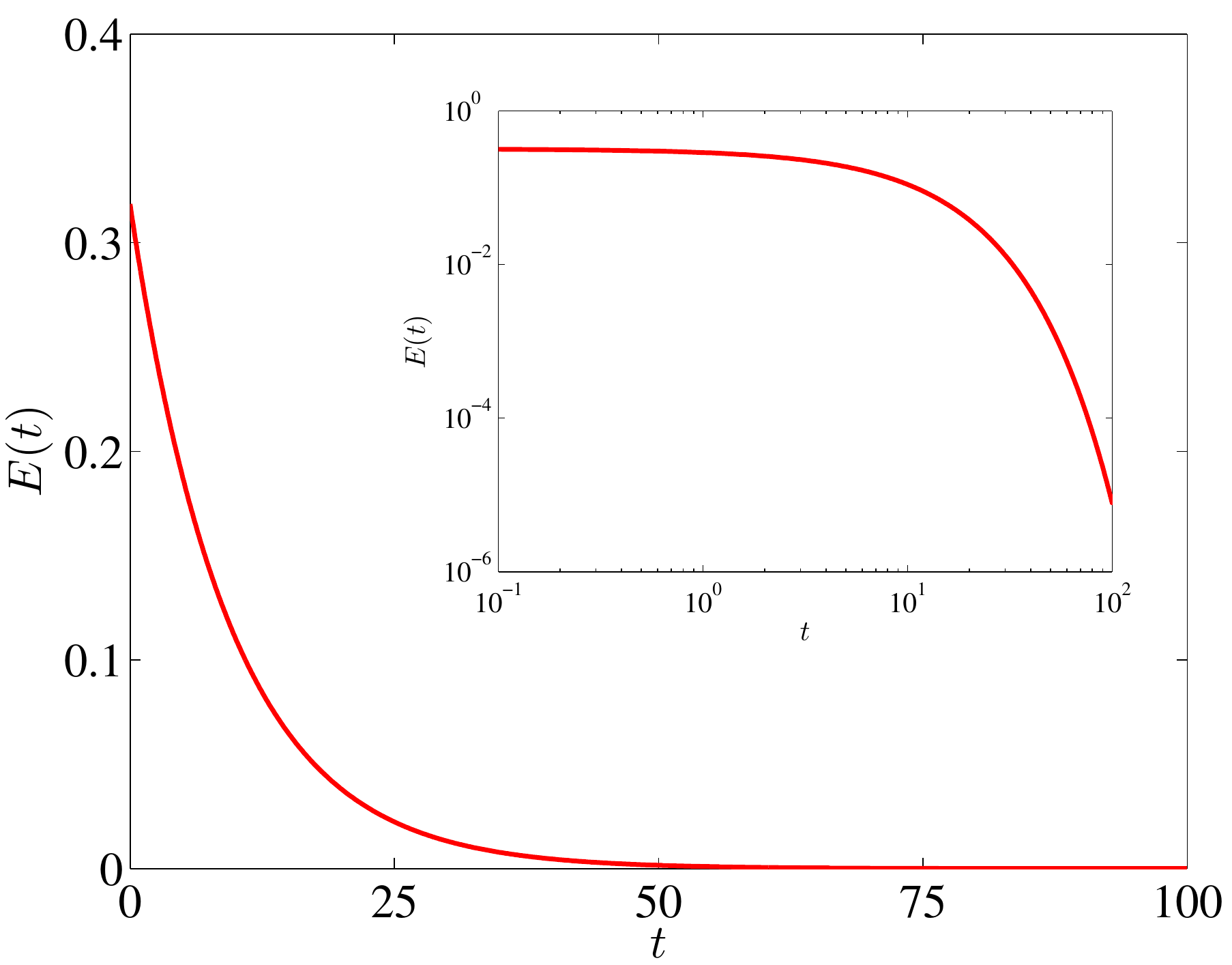}
	\end{minipage}
	\begin{minipage}{0.92\textwidth}
		\includegraphics[width=0.19\textwidth]{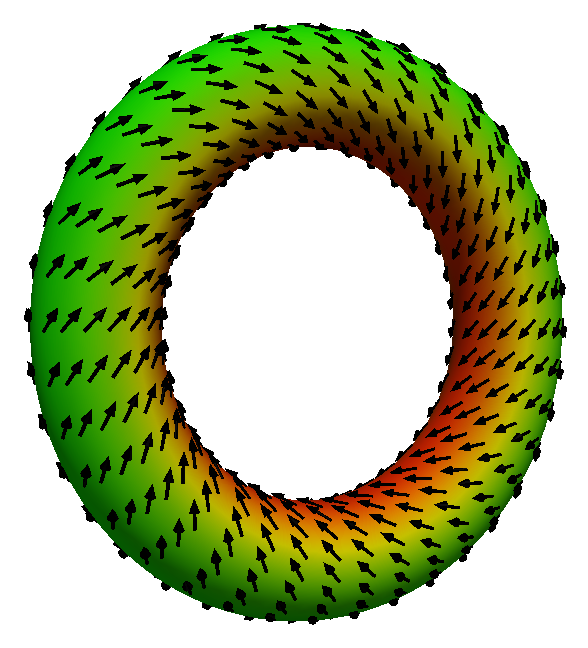}
		\includegraphics[width=0.19\textwidth]{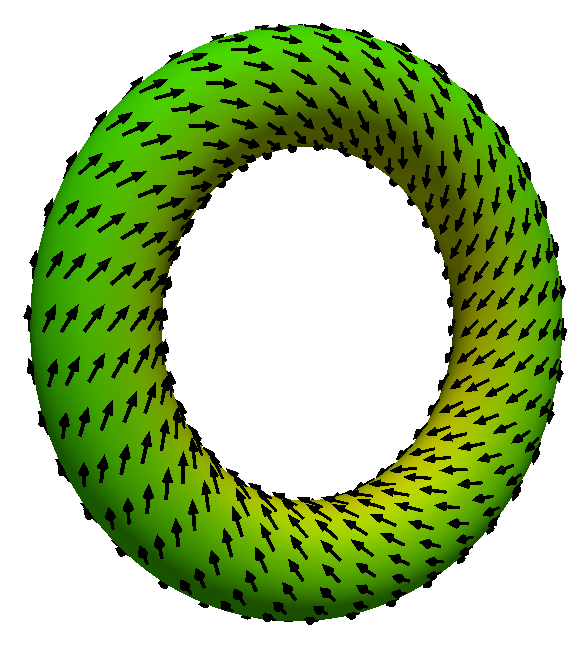}
		\includegraphics[width=0.19\textwidth]{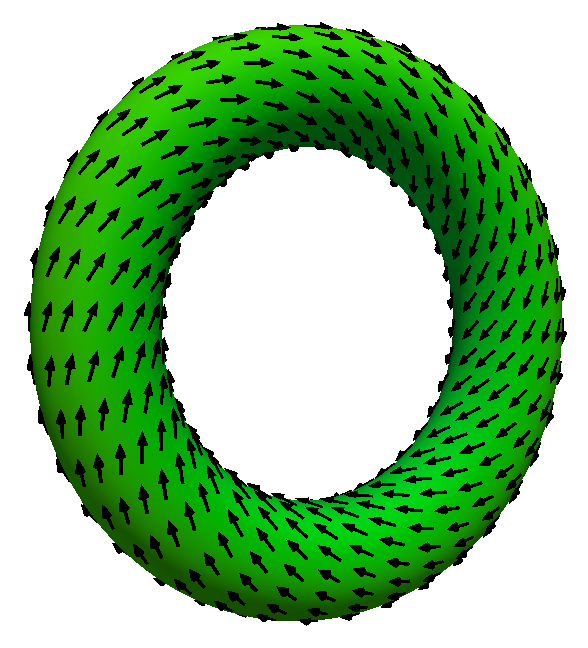}
		\includegraphics[width=0.19\textwidth]{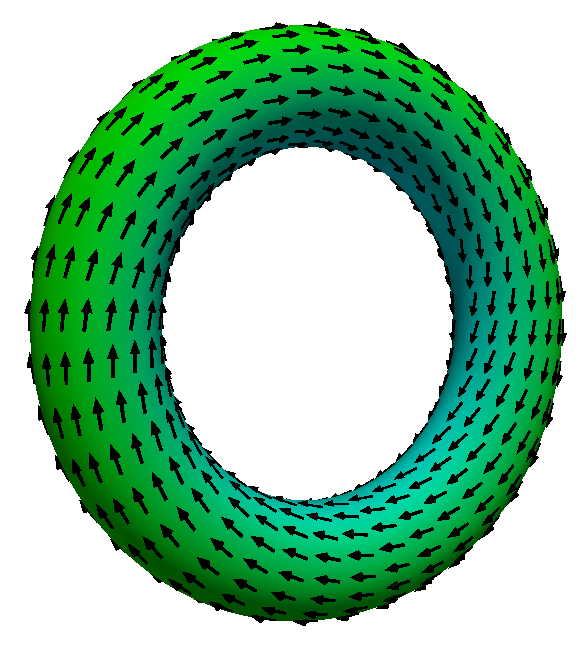}
		\includegraphics[width=0.19\textwidth]{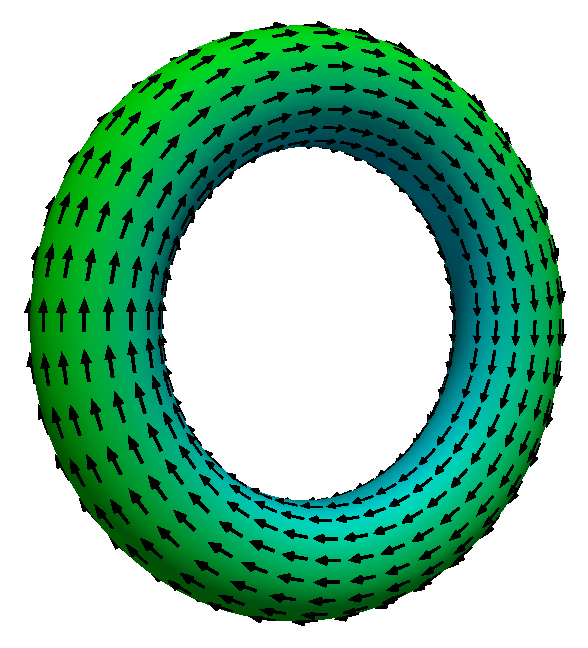}
	\end{minipage}
	\begin{minipage}{0.06\textwidth}
		\centering
		$|\vb|$\\
		\includegraphics[width=\textwidth]{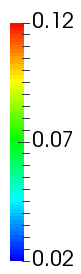}
	\end{minipage}
	\caption{(top) Harmonic vector fields $\vb_\phi^{harm}$ and $\vb_\theta^{harm}$ and kinetic energy $E$ over time $t$ for the simulation with $\vb_\theta^{harm}$ as initial condition on both equally spaced and logarithmic scales (right). (bottom) Numerical solution of $\vb$ for the simulation with $\frac{1}{2}( \vb_\phi^{harm} + \vb_\theta^{harm} )$ as initial condition computed with the DEC algorithm at $t = 0, 2, 10, 30$ and $60$ (left to right). The arrows are rescaled for better visualization.}
	\label{fig:torus:2}
\end{figure}

\subsection{Comparison}

All results for surfaces with genus $g(\surf) = 0$ demonstrate the accuracy of the DEC discretization. The plotted vortex trajectories and kinetic energy values over time are almost
indistinguishable from the SFEM results obtained by solving the vorticity-stream function formulation. The computational cost is larger for the DEC discretization, which however is also a
consequence of the implementation. Both methods are implemented in the finite element toolbox AMDiS \cite{VeyVoigt_CVS_2007,Witkowskietal_ACM_2015}, where the datastructures are optimized for SFEM, but not for DEC. A new general DEC toolbox is work in progress. 

\section{Conclusions}\label{sec5}

Even if the formulation of the incompressible surface Navier-Stokes equation is relatively old, numerical treatments on general surfaces are very rare. This also has not changed with the development of various numerical methods to solve scalar-valued partial differential equations on surfaces, such as the surface finite element method \cite{DziukElliott_IMAJNA_2007} or the diffuse interface approach \cite{Raetzetal_CMS_2006}. They are not directly applicable to vector-valued partial differential equations on surfaces. One has to define what it means 
for a vector to be parallel on the discrete representation $\SC$ of $\surf$. The concept of discrete parallel transport can be easily realized using discrete exterior calculus (DEC), see \cite{Craneetal_2013} for details. DEC thus provides an ideal framework to solve vector-valued partial differential equations on surfaces. In \cite{Nestleretal_arXiv_2016} this is shown in detail for a surface Frank-Oseen model. We here use the approach to discretize the incompressible surface Navier-Stokes equation. The discretization is based on the covariant form and utilizing a discrete version of the Hodge star $*$ and the Stokes theorem for the exterior derivative $\exd$. Non-standard in our discretization is the treatment of the discrete Hodge star and the discrete inner product. If considered in flat space the described discretization can be related to a finite difference schemes on a staggered grid. The resulting unusual stencil shows second order consistency. 

Computationally we compare results of the DEC discretization with a vorticity-stream function approach for surfaces with genus $g(\surf) = 0$. The examples use the properties of Killing vector fields and demonstrate the interplay between topology, geometry and flow properties. The numerical results are almost indistinguishable for all considered examples, varying the underlying surface $\surf$ and the Reynolds number Re. We also demonstrate the possibility to deal with harmonic vector fields using the DEC approach. It would be interesting to compare the considered vortex trajectories for larger Re with results for point vortices on closed surfaces, as e.g. considered in \cite{Dritscheletal_PRSA_2015} for ellipsoids or in \cite{Sakajoetal_PRSA_2016} for toroidal surfaces.

\begin{acknowledgement}
This work is partially supported by the German Research Foundation through grant Vo899/11. We further acknowledge computing resources provided at JSC under grant HDR06.
\end{acknowledgement}

%%%%%%%%%%%%%%%%%%%%%%%%%%%%%%%%%%%%%%%%%%%%%%%%%%%%%%%%%%%%%%%%%%%%%%%%%%%%%%%%%%%%
% Appendix sections

\section*{Appendix}
\addcontentsline{toc}{section}{Appendix}
\setcounter{section}{0}
\renewcommand\thesection{\Alph{section}}
\renewcommand*{\theHsection}{chX.\the\value{section}}

\section{Notation for DEC}
\label{app:a}

\begin{figure}[t]
  \centering
  \begin{tikzpicture}[>=latex, line width=1pt, scale=0.7]

% vertex coords
\coordinate (V1) at (-2,0);
\coordinate (V2) at (2,0);
\coordinate (V3) at (0,3);
\coordinate (V4) at (0,-3);

%points
\fill(V1) circle(2pt);
\fill(V2) circle(2pt);
\fill(V3) circle(2pt);
\fill(V4) circle(2pt);

%arrows
\draw[->] (V1) -- (V2) node [midway,above] {$e$} ;
\draw[->] (V2) -- (V3) node [midway,right] {$\tilde{e}_{2}$};
\draw[->] (V1) -- (V3) node [midway,left] {$\tilde{e}_{1}$};
\draw[->] (V4) -- (V2) node [midway,right] {$\tilde{e}_{3}$};
\draw[->] (V4) -- (V1) node [midway,left] {$\tilde{e}_{4}$};

\coordinate (VV1) at (-8,0);
\coordinate (VV2) at (-4,0);
\coordinate (VV3) at (-6,3);

\fill(VV1) circle(2pt);
\fill(VV2) circle(2pt);
\fill(VV3) circle(2pt);

\draw[->] (VV1) -- (VV2) node [midway,above] {$e$} ;
\draw[] (VV2) -- (VV3) ;
\draw[] (VV1) -- (VV3) ;

\node at (-6,1.5) {$T$};

\coordinate (VVV1) at (-8,-2);
\coordinate (VVV2) at (-4,-2);

\fill(VVV1) circle(2pt);
\fill(VVV2) circle(2pt);

\draw[->] (VVV1) -- (VVV2) node [midway,above] {$e$} ;

\node[right] at (VVV2) {$v$};

\end{tikzpicture}
	\caption{These formations always yield positive signs \( +1 \) for \( s_{\face,e} \) (top left), \( s_{v,e} \) (bottom left) and \( s_{e,\tilde{e}_{i}} \) (right) for \( i\in\left\{ 1,2,3,4 \right\}
  \), respectively. Every odd-numbered change in edge orientations results in a change of the sign \( s_{\cdot,\cdot} \).}
	\label{fig:correctEdgeOrientation}
\end{figure}
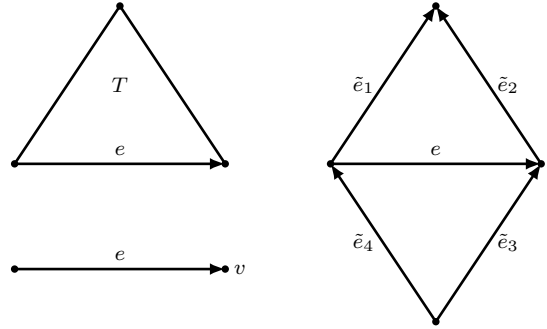

We often use the strict order relation \( \succ \) and \( \prec \) on simplices, where \( \succ \) is proverbial the ``contains'' relation, i.e. \( e \succ v \) means: the edge \( e \)
contains the vertex \( v \). Correspondingly \( \prec \) is the ``part of'' relation, i.e. \( v \prec \face \) means: the vertex \( v \) is part of the face \( \face \). Hence, we can use
this notation also for sums, like \( \sum_{f\succ e} \), i.e. the sum over all faces \( \face \) containing the edge \( e \), or \( \sum_{v\prec e} \), i.e. the sum over all vertices \( v \)
being part of the edge \( e \). Sometimes we need to determine this relation for edges more precisely with respect to the orientation. Therefore, sign functions are introduced,
\begin{align*}
	s_{\face,e} &:=
	\begin{cases}
		+1 & \text{if } e\prec \face \text{ and \( \face \) is on the left side of \( e \)} \\
		-1 & \text{if } e\prec \face \text{ and \( \face \) is on the right side of \( e \)} \formComma
	\end{cases}\\
	s_{e,\tilde{e}} &:=
	\begin{cases}
		+1 & \text{if }\measuredangle(\vect{e},\tilde{\vect{e}}) < \pi \\
		-1 & \text{if }\measuredangle(\vect{e},\tilde{\vect{e}}) > \pi
	\end{cases}\\
	s_{v,e} &:=
	\begin{cases}
		+1 & \text{if } v\prec e \text{ and \( e \) points to \( v \)} \\
		-1 & \text{if } v\prec e \text{ and \( e \) points away from \( v \)} \formComma
	\end{cases}
\end{align*}
to describe such relations between faces and edges, edges and edges or vertices and edges, respectively. Fig. \ref{fig:correctEdgeOrientation} gives a schematic illustration.

The property of a primal mesh to be well-centered ensures the existence of a Voronoi mesh (dual mesh), which is also an orientable manifold-like simplicial complex, but not well-centered.
The basis of the Voronoi mesh are not simplices, but chains of them. To identify these basic chains, we apply the (geometrical) star operator \( \star \) on the primal simplices, i.e. \(
\star v\) is the Voronoi cell corresponding to the vertex \( v \) and inherits its orientation from the orientation of the polyhedron \( \left| \SC \right| \). 
From a geometric point of view, \( \star v \) is
the convex hull of circumcenters \( c(\face) \) of all triangles \( \face\succ v \). The Voronoi edge \( \star e \) of an edge \( e \) is a connection of the right
face \( \face_{2}\succ e \) with the left face \( \face_{1}\succ e \) over the midpoint \( c(e) \). The Voronoi vertex \( \star\face \) of a face \( \face \) is simply its circumcenter \(
c(\face) \), cf. Fig.\ref{fig:VCellDecomp}. 
For a more detailed mathematical discussion see e.g. \cite{Hirani_2003,vanderzeeetal_SJSC_2010}.

The boundary operator \( \partial \) maps simplices (or chains of them) to the chain of simplices that describes its boundary with respect to its orientation (see \cite{Hirani_2003}), e.g. \( \partial(\star v)=-\sum_{e\succ v} s_{v,e}(\star e)\) (formal sum for chains) and \( \partial e = \sum_{v\prec e} s_{v,e} v \).

The expression \( \left| \cdot \right| \) measures the volume of a simplex, i.e. \( \left| \face \right| \) the area of the face \( \face \), \( \left| e \right| \) the length of the edge \( e \) and the 0-dimensional volume \( \left| v \right| \) is set to be 1. Therefore, the volume is also defined for chains and the dual mesh, since the integral is a linear functional.

\section{Second order convergence}\label{app:b}
  In this section we show that the discretization eq. \eqref{eq:lrotrotstaggered} of \( \laplaceRotRot \), defined in eq. \eqref{eq:lrotrotflat} on a staggered grid, has a truncation error of order two.
  Without loss of generality, by a quarter turn of the difference scheme in Fig. \ref{fig:laplacesStaggeredGrid} (left), we only elaborate on the discretization of
  \( (\laplaceRotRot u)^{x} \) along the horizontal x-direction.
  The first three terms in eq. \eqref{eq:lrotrotstaggered} show the well-known second order central difference approximation in vertical direction of the first term in eq. 
  \eqref{eq:lrotrotflat},
  i.e.
  \begin{align*}
    \frac{1}{h^{2}}\left( u^{x}_{i,j+1} +  u^{x}_{i,j-1} - 2u^{x}_{i,j}\right)
        &= \left( \partial_{y}^{2}u^{x} \right)^{x}_{i,j} + \landau(h^{2}) \formPeriod
  \end{align*}
  For the remaining terms, we first carry out a Taylor expansion on central vertices \( v_{i+k,j}\in\Vs \) for \( k\in\left\{ 0,1 \right\} \) in the vertical edge columns, i.e.
  \begin{align*}
      u^{y}_{i+k,j-l} 
          &= \left( u^{y} + (-1)^{l} \frac{h}{2}\partial_{y}u^{y} + \frac{h^{2}}{8}\partial_{y}^{2}u^{y} 
                    + (-1)^{l}\frac{h^{3}}{48}\partial_{y}^{3}u^{y} + \frac{h^{4}}{384}\partial_{y}^{4}u^{y}\right)_{i+k,j}\\
             &\quad +\landau(h^{5})
  \end{align*}
  for all \( l\in\left\{ 0,1 \right\} \).
  An additional horizontal expansion of sufficient order at the edge midpoint \( c(e^{x}_{i,j}) \) results in
  \begin{align*}
      \begin{aligned}
      u^{y}_{i+k,j-l}
        = &\Big(\! u^{y} + \!(-1)^{k+1}\frac{h}{2}\partial_{x}u^{y} + \frac{h^{2}}{8}\partial_{x}^{2}u^{y} 
                    + \!(-1)^{k+1}\frac{h^{3}}{48}\partial_{x}^{3}u^{y} + \frac{h^{4}}{384}\partial_{x}^{4}u^{y}  \\
                 &+ \!(-1)^{l} \frac{h}{2}\partial_{y}u^{y} + \!(-1)^{l+k+1} \frac{h^{2}}{4}\partial_{x}\partial_{y}u^{y}
                      + \!(-1)^{l} \frac{h^{3}}{16}\partial_{x}^{2}\partial_{y}u^{y} \\
                      &+ \!(-1)^{l+k+1} \frac{h^{4}}{96}\partial_{x}^{3}\partial_{y}u^{y} 
                 + \frac{h^{2}}{8}\partial_{y}^{2}u^{y} + \!(-1)^{k+1} \frac{h^{3}}{16}\partial_{x}\partial_{y}^{2}u^{y} \\
                      &+ \frac{h^{4}}{64}\partial_{x}^{2}\partial_{y}^{2}u^{y} 
                 + \!(-1)^{l}\frac{h^{3}}{48}\partial_{y}^{3}u^{y} + \!(-1)^{l+k+1} \frac{h^{4}}{96}\partial_{x}\partial_{y}^{3}u^{y} 
                 + \frac{h^{4}}{384}\partial_{y}^{4}u^{y}\!\Big)^{x}_{i,j} \\
                 &+ \landau(h^{5})
       \end{aligned}
  \end{align*}
  for all \( l,k\in\left\{ 0,1 \right\} \).
  Finally, we obtain
  \begin{align*}
      &\frac{1}{h^{2}}\left( u^{y}_{i,j} - u^{y}_{i+1,j} + u^{y}_{i+1,j-1} -u^{y}_{i,j-1} \right) \\
          &\quad= -\left( \partial_{x}\partial_{y}u^{y} 
              + \frac{h^{2}}{96}\partial_{x}\partial_{y}\left( \partial_{x}^{2}u^{y} + \partial_{y}^{2}u^{y} \right)\right)^{x}_{i,j}
              +\landau(h^{3})
  \end{align*}
  and thus a truncation error at most \( \landau(h^{2}) \) regarding \( (\laplaceRotRot u)^{x}_{i,j} \) generally.

\bibliography{biball}

\begin{thebibliography}{10}
\providecommand{\url}[1]{{#1}}
\providecommand{\urlprefix}{URL }
\expandafter\ifx\csname urlstyle\endcsname\relax
  \providecommand{\doi}[1]{DOI~\discretionary{}{}{}#1}\else
  \providecommand{\doi}{DOI~\discretionary{}{}{}\begingroup
  \urlstyle{rm}\Url}\fi

\bibitem{Marsden_1988}
Abraham, R., Marsden, J., Ratiu, T.: Manifolds, Tensor Analysis, and
  Applications.
\newblock No.~75 in Applied Mathematical Sciences. Springer (1988)

\bibitem{Arakawa_MCP_1977}
Arakawa, A., Lamb, V.: {Computational Design of the Basic Dynamical Processes
  of the UCLA General Circulation Model }.
\newblock In: General Circulation Models of the Atmosphere, pp. 173 -- 265.
  Academic Press (1977)

\bibitem{Arnoldetal_AN_2006}
Arnold, D.N., Falk, R.S., Winther, R.: Finite element exterior calculus,
  homological techniques, and applications.
\newblock Acta Numerica \textbf{15}, 1--155 (2006)

\bibitem{ArroyoDeSimone_PRE_2009}
Arroyo, M., DeSimone, A.: Relaxation dynamics of fluid membranes.
\newblock Physical Review E \textbf{79}, 031,915 (2009)

\bibitem{Barrettetal_PRE_2015}
Barrett, J., Garcke, H., N\"urnberg, R.: Numerical computations of the dynamics
  of fluidic membranes and vesicles.
\newblock Physical Review E \textbf{92}, 052,704 (2015)

\bibitem{Botheetal_JFM_2010}
Bothe, D., Pr{\"u}ss, J.: On the two-phase {Navier-Stokes} equations with
  {Boussinesq-Scriven} surface.
\newblock Journal of Mathematical Fluid Mechanics \textbf{12}, 133--150 (2010)

\bibitem{Craneetal_2013}
Crane, K., de~Goes, F., Desbrun, M., Schr{\"o}der, P.: {Digital geometry
  processing with discrete exterior calculus}.
\newblock In: ACM SIGGRAPH Courses, pp. 1--126 (2013)

\bibitem{Desbrunetal_arXiv_2005}
Desbrun, M., Hirani, A., Leok, M., Marsden, J.: Discrete exterior calculus.
\newblock arXiv:math/0508341  (2005)

\bibitem{Doerriesetal_PRE_1996}
D\"orries, G., Foltin, G.: Energy dissipation of fluid membranes.
\newblock Physical Review E \textbf{53}, 2547--2550 (1996)

\bibitem{Dritscheletal_PRSA_2015}
Dritschel, D.G., Boatto, S.: The motionof point vortices on closed surfaces.
\newblock Proceedings of the Royal Society A \textbf{471}, 20140,890 (2015)

\bibitem{DziukElliott_JCM_2007}
Dziuk, G., Elliott, C.: Surface finite elements for parabolic equations.
\newblock Journal of Computational Mathematics \textbf{25}, 385--407 (2007)

\bibitem{DziukElliott_IMAJNA_2007}
Dziuk, G., Elliott, C.M.: Finite elements on evolving surfaces.
\newblock IMA Journal of Numerical Analysis \textbf{27}, 262--292 (2007)

\bibitem{EbinMarsden_AM_1970}
Ebin, D.G., Marsden, J.: Groups of diffeomorphisms and the motion of an
  incompressible fluid.
\newblock Annals of Mathematics \textbf{92}, 102--163 (1970)

\bibitem{Elcottetal_2007}
Elcott, S., Tong, Y., Kanso, E., Schr{\"o}der, P., Desbrun, M.: Stable,
  circulation-preserving, simplicial fluids.
\newblock ACM Transactions on Graphics \textbf{26}, 4 (2007)

\bibitem{Fanetal_JCP_2010}
Fan, J., Han, T., Haataja, M.: Hydrodynamic effects on spinodal decomposition
  kinetics in planar lipid bilayer membranes.
\newblock The Journal of Chemical Physics \textbf{133}, 235,101 (2010)

\bibitem{Fisheretal_2006}
Fisher, M., Springborn, B., Bobenko, A., Schr{\"o}der, P.: {An algorithm for
  the construction of intrinsic Delaunay triangulations with applications to
  digital geometry processing}.
\newblock In: ACM SIGGRAPH Courses, pp. 69--74 (2006)

\bibitem{Gortleretal_2006}
Gortler, S., Gotsman, C., Thurston, D.: Discrete one-forms on meshes and
  applications to {3D} mesh parameterization.
\newblock Computer Aided Geometric Design \textbf{33}, 83--112 (2006)

\bibitem{Griebeletal_SPP_2017}
Griebel, M., Rieger, C., Schier, A.: ?
\newblock ? \textbf{?}, ? (2017)

\bibitem{Guetal_2003}
Gu, X., Yau, S.T.: {Global conformal surface parameterization}.
\newblock In: ACM/EG Symposium on Geometry Processing, pp. 127--137 (2003)

\bibitem{Hirani_2003}
Hirani, A.N.: Discrete exterior calculus.
\newblock Ph.D. thesis, California Institute of Technology, Pasadena, CA, USA
  (2003)

\bibitem{Huetal_PRE_2007}
Hu, D., Zhang, P., E, W.: Continuum theory of a moving membrane.
\newblock Physical Review E \textbf{75}, 041,605 (2007)

\bibitem{Mercat_CMP_2001}
Mercat, C.: Discrete {R}iemann surfaces and the {I}sing model.
\newblock Communications in Mathematical Physics \textbf{218}, 177--216 (2001)

\bibitem{MitreaTaylor_MA_2001}
Mitrea, M., Taylor, M.: Navier-{S}tokes equations on {Lipschitz} domains in
  {Riemannian} manifolds.
\newblock Mathematische Annalen \textbf{321}, 955--987 (2001)

\bibitem{Mohamed_CAD_2016}
Mohamed, M.S., Hirani, A.N., Samtaney, R.: Comparison of discrete hodge star
  operators for surfaces.
\newblock Computer-Aided Design  (2016).
\newblock \doi{10.1016/j.cad.2016.05.002}

\bibitem{Mohamedetal_JCP_2016}
Mohamed, M.S., Hirani, A.N., Samtaney, R.: Discrete exterior calculus
  discretization of incompressible {Navier-€"Stokes} equations over surface
  simplicial meshes.
\newblock Journal of Computational Physics \textbf{312}, 175 -- 191 (2016)

\bibitem{Mullenetal_2009}
Mullen, P., Crane, K., Pavlov, D., Tong, Y., Desbrun, M.: Energy-preserving
  integrators for fluid animation.
\newblock ACM Transactions on Graphics \textbf{28}, 38 (2009)

\bibitem{Nestleretal_arXiv_2016}
{Nestler}, M., {Nitschke}, I., {Praetorius}, S., {Voigt}, A.: {Orientational
  order on surfaces - the coupling of topology, geometry and dynamics}.
\newblock arXiv:1608.01343  (2016)

\bibitem{Nitschkeetal_inprep}
Nitschke, I., Voigt, A.: Curvature approximation of discrete surfaces - a
  discrete exterior calculus approach.
\newblock in preparation

\bibitem{Nitschkeetal_JFM_2012}
Nitschke, I., Voigt, A., Wensch, J.: A finite element approach to
  incompressible two-phase flow on manifolds.
\newblock Journal of Fluid Mechanics \textbf{708}, 418--438 (2012)

\bibitem{Polthieretal_2003}
Polthier, K., Preu{\ss}, E.: {Identifying vector field singularities using a
  discrete Hodge decomposition}.
\newblock In: H.~Hege, K.~Polthier (eds.) Visualization and Mathematics III,
  pp. 113--134. Springer (2003)

\bibitem{Raetzetal_CMS_2006}
R\"{a}tz, A., Voigt, A.: {PDE's} on surfaces: {A} diffuse interface approach.
\newblock Communications in Mathematical Sciences \textbf{4}, 575--590 (2006)

\bibitem{Reutheretal_MMS_2015}
Reuther, S., Voigt, A.: The interplay of curvature and vortices in flow on
  curved surfaces.
\newblock Multiscale Modeling \& Simulation \textbf{13}, 632--643 (2015)

\bibitem{Reutheretal_JCP_2016}
Reuther, S., Voigt, A.: {Incompressible two-phase flows with an inextensible
  Newtonian fluid interface}.
\newblock Journal of Computational Physics \textbf{322}, 850--858 (2016)

\bibitem{Sakajoetal_PRSA_2016}
Sakajo, T., Shimizu, Y.: Point vortex interactions on a toroidal surface.
\newblock Proceedings of the Royal Society A \textbf{472}, 20160,271 (2016)

\bibitem{Scriven_CES_1960}
Scriven, L.E.: Dynamics of a fluid interface equation of motion for {Newtonian}
  surface fluids.
\newblock Chemical Engineering Science \textbf{12}, 98--108 (1960)

\bibitem{Secombetal_QJMAM_1982}
Secomb, T.W., Skalak, R.: Surface flow of viscoelastic membranes in viscous
  fluids.
\newblock The Quarterly Journal of Mechanics and Applied Mathematics
  \textbf{35}, 233--247 (1982)

\bibitem{Tongetal_2006}
Tong, Y., Alliez, P., Cohen-Steiner, D., Desbrun, M.: {Designing
  quadrangulations with discrete harmonic forms}.
\newblock In: ACM/EG Symposium on Geometry Processing, pp. 201--210 (2006)

\bibitem{Tongetal_2003}
Tong, Y., Lombeyda, S., Hirani, A.N., Desbrun, M.: Discrete multiscale vector
  field decomposition.
\newblock ACM Transactions on Graphics \textbf{22}, 445--452 (2003)

\bibitem{vanderzeeetal_SJSC_2010}
VanderZee, E., Hirani, A.N., Guoy, D., Ramos, E.A.: Well-centered
  triangulation.
\newblock SIAM Journal on Scientific Computing \textbf{31}, 4497--4523 (2010)

\bibitem{Vaxmanetal_2016}
Vaxman, A., Campen, M., Diamanti, O., Panozzo, D., Bommes, D., Hildebrandt, K.,
  Ben-Chen, M.: {Directional field synthesis, design and processing}.
\newblock In: EUROGRAPHICS - STAR, vol.~35, pp. 1--28 (2016)

\bibitem{VeyVoigt_CVS_2007}
Vey, S., Voigt, A.: {AMDiS: Adaptive multidimensional simulations}.
\newblock Computing and Visualization in Science \textbf{10}, 57--67 (2007)

\bibitem{Witkowskietal_ACM_2015}
Witkowski, T., Ling, S., Praetorius, S., Voigt, A.: Software concepts and
  numerical algorithms for a scalable adaptive parallel finite element method.
\newblock Advances in Computational Mathematics \textbf{41}, 1145--1177 (2015)

\end{thebibliography}
\bibliographystyle{spmpsci}
% or without title 
%\bibliographystyle{spphys}

\end{document}